\documentclass[a4paper, 11pt, twoside]{article}
\usepackage{fullpage}

\usepackage{amsmath,amsfonts,amssymb,amsthm,array}
\usepackage{algorithm}
\usepackage{caption}
\usepackage{subcaption}
\usepackage{algpseudocode}
\usepackage{bm}
\usepackage{color}
\usepackage{graphicx}
\usepackage{enumitem}
\usepackage{nicefrac}

\usepackage[utf8]{inputenc} % allow utf-8 input
\usepackage[T1]{fontenc}    % use 8-bit T1 fonts
\usepackage{hyperref}       % hyperlinks
\usepackage{url}            % simple URL typesetting
\usepackage{booktabs}       % professional-quality tables
\usepackage{amsfonts}       % blackboard math symbols
\usepackage{nicefrac}       % compact symbols for 1/2, etc.
\usepackage{microtype}      % microtypography
\usepackage{xcolor}         % colors
\usepackage{amsmath,amsfonts,amssymb,amsthm,array}
\usepackage{comment}
\usepackage{algorithm}
\usepackage[font=small]{caption}
\usepackage{subcaption}
\usepackage{algpseudocode}
\usepackage{bm}
\usepackage{color}
\usepackage{graphicx}
\usepackage{enumitem}
\usepackage{xargs}
\usepackage{multirow}
\usepackage{makecell}
\usepackage{wrapfig}
\usepackage{bbm, dsfont}
\usepackage[colorinlistoftodos,bordercolor=orange,backgroundcolor=orange!20,linecolor=orange,textsize=scriptsize]{todonotes}

\usepackage{xspace}

\usepackage{cleveref}

\usepackage{makecell}
\usepackage{caption}
\usepackage{multirow}

\usepackage{colortbl}
\definecolor{bgcolor}{rgb}{0.8,1,1}
\definecolor{bgcolor2}{rgb}{0.8,1,0.8}
\usepackage{threeparttable}

\usepackage{tcolorbox}
\usepackage{pifont}
\definecolor{mydarkgreen}{RGB}{39,130,67}
\definecolor{mydarkred}{RGB}{192,47,25}

\usepackage{subcaption}
\usepackage{wrapfig}

\allowdisplaybreaks
\newcommand{\R}{\ensuremath{\mathbb{R}}}

\newcommand{\EE}{\mathbb{E}}

\def\<#1,#2>{\left\langle #1,#2 \right\rangle}
\newcommand{\EndProof}{\begin{flushright}$\square$\end{flushright}}

\newcommandx{\norm}[2][2=]{\Vert#1 \Vert_{{#2}}}

\newcommand{\dobrush}{\mathsf{\Delta}}
\newcommandx{\dobru}[3][1=,3=]{\dobrush_{#1}^{#3}( #2)}  %%% dobrushin coefficient

\def\L{{\cal L}}
\def\R{\mathbb{R}}

\def\R{\mathbb R}

\def\EE{\mathbb E}

\def\la{\langle}
\def\ra{\rangle}

\usepackage{caption}
\usepackage{xcolor}
\definecolor{ForestGreen}{RGB}{34,139,34}

\allowdisplaybreaks
\usepackage{mathtools}
\newcommand{\eqdef}{\vcentcolon=}

\usepackage{colortbl}
\definecolor{bgcolor2}{rgb}{0.8,1,1}
\definecolor{bgcolor}{rgb}{0.8,1,0.8}
\definecolor{niceblue}{rgb}{0.0,0.19,0.56}

\def\L{{\cal L}}

\usepackage{tikz}

\theoremstyle{plain}

\theoremstyle{remark}
\usepackage{mdframed} 
\usepackage{thmtools}

\usepackage{caption}
\usepackage{multirow}
\usepackage{colortbl}
\definecolor{bgcolor2}{rgb}{0.8,1,1}
\definecolor{bgcolor}{rgb}{0.8,1,0.8}

\usepackage{hyperref}
\hypersetup{colorlinks,linkcolor={blue},citecolor={niceblue},urlcolor={blue}}

\usepackage{natbib}
\usepackage{sidecap}

\definecolor{shadecolor}{gray}{0.9}
\declaretheoremstyle[
headfont=\normalfont\bfseries,
notefont=\mdseries, notebraces={(}{)},
bodyfont=\normalfont,
postheadspace=0.5em,
spaceabove=1pt,
mdframed={
  skipabove=8pt,
  skipbelow=8pt,
  hidealllines=true,
  backgroundcolor={shadecolor},
  innerleftmargin=4pt,
  innerrightmargin=4pt}
]{shaded}

\declaretheorem[style=shaded,within=section]{definition}
\declaretheorem[style=shaded]{theorem}

\declaretheorem[name=A, style=shaded]{assumption}
\declaretheorem[style=shaded]{corollary}

\declaretheorem[style=shaded]{lemma}

\usepackage{xspace}

\usepackage[colorinlistoftodos,bordercolor=orange,backgroundcolor=orange!20,linecolor=orange,textsize=scriptsize]{todonotes}

\usepackage{amsmath}
\usepackage{amssymb}
\usepackage{mathtools}
\usepackage{amsthm}

\usepackage{verbatim}
\usepackage[utf8]{inputenc} % allow utf-8 input
\usepackage[T1]{fontenc}    % use 8-bit T1 fonts
\usepackage{hyperref}       % hyperlinks
\usepackage{url}            % simple URL typesetting
\usepackage{booktabs}       % professional-quality tables
\usepackage{amsfonts}       % blackboard math symbols
\usepackage{nicefrac}       % compact symbols for 1/2, etc.
\usepackage{microtype}      % microtypography
\usepackage{xcolor}         % colors
\usepackage{amsmath}
\usepackage{amssymb}
\usepackage[colorinlistoftodos,bordercolor=orange,backgroundcolor=orange!20,linecolor=orange,textsize=scriptsize]{todonotes}
\usepackage{enumitem}
\usepackage{pifont}
\usepackage{algorithm,algorithmicx,algpseudocode}
\usepackage{subcaption}
\usepackage{wrapfig}

\usepackage{makecell}
\usepackage{caption}
\usepackage{multirow}

\title{Accelerated Methods with Compression for Horizontal and Vertical Federated Learning}

\begin{document}

\title{\textbf{Accelerated Methods with Compression for Horizontal and Vertical Federated Learning}
}

\author{Sergey Stanko $^{1}$
\quad Timur Karimullin $^{1}$
\quad Aleksandr Beznosikov $^{1,2}$
\quad Alexander Gasnikov $^{1}$
}
\date{$^1$ Moscow Institute of Physics and Technology, Russian Federation \\
$^2$ Mohamed bin Zayed University of Artificial Intelligence, United Arab Emirates
}
\maketitle
\begin{abstract}
Distributed optimization algorithms have emerged as a superior approaches for solving machine learning problems. To accommodate the diverse ways in which data can be stored across devices, these methods must be adaptable to a wide range of situations. As a result, two orthogonal regimes of distributed algorithms are distinguished: horizontal and vertical. During parallel training, communication between nodes can become a critical bottleneck, particularly for high-dimensional and over-parameterized models. Therefore, it is crucial to enhance current methods with strategies that minimize the amount of data transmitted during training while still achieving a model of similar quality. This paper introduces two accelerated algorithms with various compressors, working in the regime of horizontal and vertical data division. By utilizing a momentum and variance reduction technique from the \texttt{Katyusha} algorithm, we were able to achieve acceleration and demonstrate one of the best asymptotics for the horizontal case. Additionally, we provide one of the first theoretical convergence guarantees for the vertical regime. Our experiments involved several compressor operators, including RandK and PermK, and we were able to demonstrate superior practical performance compared to other popular approaches.
\end{abstract}

\section{Introduction}\label{sec:intro}
As machine learning continues to gain popularity, it is essential for training algorithms to maintain exceptionally high performance in order to effectively address real-world issues. To increase the quality of models' predictions, larger datasets are used. However, training models on a huge amount of samples on a single machine can be very time-consuming, which is why faster optimization approaches are needed.

To address this problem, the scientific community came to distributed algorithms, where the calculation process is divided among different devices (also called machines, nodes, clusters, or workers). Such parallel computation can be used in situations, where data is distributed across several machines, as in the cases of classical distributed \cite{verbraeken2020survey} and federated learning approaches \citep{konečný2017federated}. The latter can be divided into two different regimes. One of these is horizontal federated learning, which is the most popular and extensively studied approach. In this scenario, each worker possesses their own collection of samples, but share the same set of features. A different situation is considered in the vertical case \citep{Zhang_2021}, where every device has a unique set of features of the same samples. The practical use of the latter regime can be seen in a simple example, where a digital finance company wants to evaluate the risk of approving a loan. To make a qualified prediction, it needs to collect different types of data (features) from the same person (sample). For example, online shopping information from an E-commerce company or average monthly deposit and account information from a bank. However, legal restrictions or competition between participants prevent these organizations from sharing data, which is why a vertical case is suitable in this situation \citep{gu2020privacypreserving}. But classical distributed and horizontal federative learning methods can not be forgotten, as they can show a great performance in accelerating training speed \citep{goyal2018accurate}. This is why distributed optimization is gaining a lot of attention nowadays.

However, parallelization of a task practically does not lead to an ideal decrease in time. That means that having $N$ devices does not accelerate task by $N$ times. This happens because of limited ability of networks to exchange information. Thus, the key bottleneck of parallel computation is the communication part. There have been considered several ways of dealing with this issue \cite{konevcny2016federated, smith2018cocoa}, but in our paper we concentrate solely on reducing communication cost of each iteration by decreasing the size of sending information also known as compression technique \cite{chilimbi2014project, alistarh2017qsgd}. Due to reducing quality of communicated information, gradient descent methods perform greater amount of iterations, but overall time complexity can be reduced. 

Compression mathematically can be represented in the form of a vector function, which is called compression operator. It takes a vector and returns a vector of the same dimension, but the last usually takes less resources in communication networks. In our paper, we consider such operators with the following properties:

\begin{definition} \label{ass:ass000}
We say that the compression operator $Q$ is unbiased, if
\begin{equation*}
\mathbb{E} \left [Q(x) \right ] = x \text{ , $\forall x \in \mathbb{R}^d$}.
\end{equation*}

We also assume that there exists a constant $w \geq 1$, such:
\begin{equation*}
\mathbb{E} \left [\| Q(x)\|^2 \right ] \leq \omega \| x \|^2 \text{ , $\forall x \in \mathbb{R}^d$}.
\end{equation*}

We say that the compressor operator $Q$ has density coefficient $\beta$, if $Q(x)$ requires in $\beta$ times less space complexity than $x$. 
\end{definition}

\renewcommand{\arraystretch}{2}
\begin{table*}[ht]
\vspace{-0.3cm}
    \centering
    \small
%    \scriptsize
\captionof{table}{Summary of bounds for iteration complexities for finding an $\varepsilon$-solution. Convergence is measured by the distance to the solution.}
\vspace{-0.2cm}
    \label{tab:comparison0}   
    \scriptsize
    \resizebox{\linewidth}{!}{
  \begin{threeparttable}
    \begin{tabular}{|c|c|c|}
    \hline
     \textbf{Regime} & \textbf{Reference } & \textbf{Iteration complexity}\\ 
     \cline{1-3}
     \multirow{3}{*}{\rotatebox[origin=c]{90}{\textbf{Horizontal \quad\quad\quad\quad\quad}}} & \texttt{QSGD} \cite{alistarh2017qsgd} \tnote{{\color{blue}(1)}} & $\mathcal{O}\left(\frac{L}{\mu}\left(1 + \frac{\omega - 1}{n} \right)\log{\frac{1}{\varepsilon}}\right)$ \label{que:qsgd}\\
     \cline{2-3}
     & \texttt{DIANA} \cite{mishchenko2023distributed} & $\mathcal{O}\left(\left(\frac{L}{\mu}\left(1 + \frac{2(\omega - 1)}{n}\right) + \omega  \right)\log{\frac{1}{\varepsilon}}\right)$ \\
     \cline{2-3}
     & \texttt{VR-DIANA} \cite{horváth2019stochastic} & $\mathcal{O}\left(\left(\frac{L}{\mu}\left(1 + \frac{\omega - 1}{n}\right) + \omega  \right)\log{\frac{1}{\varepsilon}}\right)$\\
     \cline{2-3}
     & \texttt{ADIANA} \cite{adiana} \tnote{{\color{blue}(2)}} & $\mathcal{O}\left(\left( \sqrt{\frac{L}{\mu}}\left(1 + \sqrt{\left ( \frac{\omega - 1}{n} + \sqrt{\frac{\omega - 1}{n}}  \right) \omega} \right) + \omega\right)\log{\frac{1}{\varepsilon}}\right)$ \\
     \cline{2-3}
     & \texttt{ADIANA} \cite{he2024unbiased} \tnote{{\color{blue}(3)}}  & $\mathcal{O}\left(\left(\sqrt{\frac{L}{\mu}}\left(\sqrt{\frac{\omega^2}{n}} + 1\right) + \omega\right) \log{\frac{1}{\varepsilon}}\right)$ \\
     \cline{2-3}
     & \texttt{MARINA} \cite{marina} \tnote{{\color{blue}(4)}} & $\mathcal{O}\left(\left(\frac{L}{\mu}\left(1 + \frac{\omega - 1}{\sqrt{n}}\right) + \omega\right)\log{\frac{1}{\varepsilon}}\right)$\\
     \cline{2-3} 
     & \texttt{MASHA} \cite{beznosikov2023distributed} \tnote{{\color{blue}(5)}} & $\mathcal{O}\left(\left(\frac{L}{\mu}\sqrt{\left(w + \frac{(\omega - 1)^2}{n}\right)} + \omega  \right)\log{\frac{1}{\varepsilon}}\right)$\\
     \cline{2-3}
     & \texttt{Three Pillars Algorithm} \cite{beznosikov2023similarity} \tnote{{\color{blue}(5, 6)}}  &  $\mathcal{O}\left(\left(\frac{L}{\mu}\sqrt{n} + n\right)\log{\frac{1}{\varepsilon}}\right)$\\
     \cline{2-3}
     &\label{tab:one} \cellcolor{bgcolor2}{Algorithm \ref{alg:one}(This paper)} & \cellcolor{bgcolor2}{$\mathcal{O}\left(\left(\sqrt{\frac{L}{\mu}}\left(\sqrt{\frac{\omega^2}{n}} + 1\right) + \omega\right) \log{\frac{1}{\varepsilon}}\right)$} \\
     \hline
     \multirow{2}{*}{\rotatebox[origin=c]{90}{\textbf{Vertical \quad\quad}}} & \texttt{AVFL} \cite{Cai_2022}; \texttt{CE-VFL} \cite{sun2023communicationefficient};  & \multirow{2}{*}{No theoretical results}\\
     & \texttt{SecureBoost+} \cite{chen2021secureboost}; \texttt{eHE-SecureBoost} \cite{xu2021efficient}; & \\
     \cline{2-3}
     & \texttt{CVFL} \cite{castiglia2023compressedvfl}  & No concrete number of iterations \tnote{{\color{blue}(7)}}\\
     \cline{2-3}
     & \cellcolor{bgcolor2}{Algorithm \ref{alg:orig}(This paper)}& \cellcolor{bgcolor2}{$\mathcal{O}\left(\left(\sqrt{\frac{\bar{L}}{\mu}s} + s \right) \log{\frac{1}{\varepsilon}}\right)$}\\
     \hline
     
    \end{tabular}
   \begin{tablenotes}
    {\small
    \item [] \tnote{{\color{blue}(1)}} correct tuning of step size, no convergence with fixed step size;
    \tnote{{\color{blue}(2)}} this is the complexity derived in the original paper \cite{adiana};
    \tnote{{\color{blue}(3)}} this is the complexity derived by a refined analysis in the preprint \cite{he2024unbiased}
    \tnote{{\color{blue}(4)}} under P\L \ condition;
    \tnote{{\color{blue}(5)}} for VI and SPPs;
    \tnote{{\color{blue}(6)}} for PermK compressor;
    \tnote{{\color{blue}(7)}} for special compressor convergence rate is $\mathcal{O}(\frac{1}{\sqrt{K}})$, no guaranties in the case of arbitrary compressor;
     \item [] {\em Notation:} $\mu$ = constant of strong convexity, $L$ = smoothness constant of the target function, $\omega$ = compression constant (see Definition \ref{ass:ass000}), $n$ =  number of workers, $s$ = number of samples, $\bar{L} = \frac{1}{s}\sum\limits_{i = 1}^sL_i$
    }
\end{tablenotes}    
    \end{threeparttable}
    }
\end{table*}

\vspace{-0.2cm}
\subsection{Related works}
\vspace{-0.2cm}
\textbf{Unbiased compression.} An idea to shrink the vector was researched in \citep{doi:10.1137/100802001}. This work developed a variant of single-node gradient descent in which only some random coordinates are updated at every step. Later, in the literature on compression, operators of this kind began to be called “RandK” \citep{beznosikov2024biased}. In the paper by \citep{richtárik2013distributed}, coordinate descent method was adapted for distributed optimization. Generalization with arbitrary unbiased compressor was firstly introduced in \texttt{QSGD} paper \citep{alistarh2017qsgd}, but this method does not converge to the true solution, but rather to some neighbourhood, dependent on functions' variance \citep{gorbunov2021distributed}.

\textbf{Variance reduction.} The problem mentioned above is characterized by a non-zero limit in the optimal point of the variance of stochastic gradient. A technique of choosing a vector of step direction to fix this challenge is called variance reduction. Firstly used in standard optimization problem \cite{NIPS2013_ac1dd209, nguyen2017sarah} and later was adapted for distributed optimization in \texttt{DIANA} \citep{mishchenko2023distributed}, by compressing not gradient itself, but rather gradient difference, and later for non-convex distributed case in \texttt{MARINA} \citep{marina}. This technique is also used in variational inequalities, for example, in \citep{beznosikov2023distributed}.

\textbf{Acceleration.} The optimal algorithm for a non-distributed strongly convex problem, which uses variance reduction technique, is \texttt{Katyusha} \citep{BadKatyusha}. It is based on acceleration technique, called by authors “negative momentum” with the combination of stochastic adjustment of the gradient in the old point. Vector, in which gradient is calculated, updates with the loop, which is theoretically no more effective than an update with a certain probability represented in \texttt{L-Katyusha} \citep{kovalev2019dont}, but is difficult to perceive and empirically slower.

\vspace{-0.2cm}
\subsection{Our contributions}
\vspace{-0.2cm}
$\bullet$ \textbf{New distributed algorithms.} We present new distributed algorithms with compression for various regimes of data division. As mentioned before, parallel algorithms can be based on non-parallel approaches with variance reduction, e.g.: for non-convex problem, \texttt{MARINA} \citep{marina} is based on \texttt{PAGE} \citep{li2021page}, or for variational inequalities, \texttt{MASHA} \cite{beznosikov2023distributed} is based on \citep{alacaoglu2022stochastic}. In our paper, we base our results on L-Katyusha \citep{kovalev2019dont}, one of the state-of-the-art algorithm for strongly convex problem with variance reduction technique.

$\bullet$ \textbf{Horizontal regime.} We propose \texttt{L-Katyusha}-based algorithm with horizontal data division -- \texttt{DHPL-Katyusha} (\texttt{Distributed Horizontally Partitioned L-Katyusha} from Algorithm \ref{alg:one}), in which every worker compresses and sends the difference between gradients of its own function in new and old points, after which regular \texttt{L-Katyusha} algorithm is performed locally. The asymptotics of our algorithm is compared with other popular approaches in Table \ref{tab:comparison0}. Unlike \texttt{QSGD} \citep{alistarh2017qsgd}, \texttt{DIANA} \citep{mishchenko2023distributed}, \texttt{VR-DIANA} \citep{horváth2019stochastic}, \texttt{MASHA} \citep{beznosikov2023distributed} and  \texttt{Three Pillars Algorithm} \citep{beznosikov2023similarity}, our algorithm is accelerated (e.g. has asymptotics proportional to $\sqrt{\frac{L}{\mu}}$, where $L$ is a smoothness constant and $\mu$ is a constant of strong convexity), and in comparison with the improved version of \texttt{ADIANA}, \texttt{DHPL-Katyusha} has the same asymptotics, but shows better empirical results in experiments.

$\bullet$ \textbf{Vertical regime.} We strive to present \texttt{DVPL-Katyusha} and \texttt{DVPL-Katyusha with scalar compression} (Algorithm \ref{alg:orig} and \ref{alg:two} respectively) — algorithms with compression working under the assumption of the vertical data division. The first one uses RandK compressor, and the second utilize MSE loss and arbitrary (with respect to Definition \ref{ass:ass000}) compressor. In both algorithms every worker sends its own part of scalar values, needed to calculate total loss, and then perform other part of \texttt{L-Katyusha} locally with its own set of features. We estimate the asymptotics in such settings and show a lower bound for scalar compression on a constructive example. We are one of the first to present algorithm with vertical regime and guaranties for theoretical convergence. Other approaches, represented in Table \ref{tab:comparison0}, such as \texttt{AVFL} \citep{Cai_2022}, \texttt{CE-VFL} \citep{sun2023communicationefficient}, \texttt{SecureBoost+} \citep{chen2021secureboost} and \texttt{eHE-SecureBoost} \citep{xu2021efficient} do not have it. The \texttt{CVFL} algorithm \citep{castiglia2023compressedvfl} has theoretical convergence rates only for special compressors, but also lacks concrete convergence guarantees.

$\bullet$ \textbf{Various compressors.} We obtain similar estimates for different compressors. By obtaining $\omega$ from Definition \ref{ass:ass000}, we research RandK operator, which retains only $K$ random coordinates from a vector. In this paper, we also implement one of the state-of-the-art compressor PermK \citep{szlendak2021permutation}, which distributes random permutations of vector components across all workers to compute only on them. As PermK can be represented in the form of several correlated compressors, $\omega$ can not be defined for it. Thus, a slightly different analysis is used in order to develop a proof.

$\bullet$ \textbf{Numerical experiments.} In numerical experiments, we illustrate the most important properties of the new methods. The results correspond to the theory developed.

\vspace{-0.2cm}
\subsection{Technical preliminaries}
\vspace{-0.2cm}

\textbf{Notations.} $\la x,y \ra \eqdef \sum_{i=1}^d x_i y_i$ is used to denote standard inner product of $x,y\in\R^d$. Operator $\EE[\cdot]$ denotes mathematical expectation. We denote $A^T_i \in \R^d$ as the $i$-th column of the matrix $A^T$.

\vspace{-0.2cm}
\section{Horizontal Case }
\vspace{-0.2cm}
In this section, we provide an algorithm for the horizontal division of data, which we call \texttt{DHPL-Katyusha} (Algorithm \ref{alg:one}). We state the standard distributed learning problem: each worker has its own dataset, on which it can calculate the loss. The target loss function, which we need to minimize, is obtained by summarizing and averaging all workers' losses. Formally, it can be written in the form of: 
\begin{equation} \label{eq:hor}
\min_{x \in \mathbb{R}^d} \left[ f(x) \eqdef \frac{1}{n} \sum\limits_{i = 1}^n f_i(x) \right].
\end{equation}

Communications between nodes can be centralized, e.g. they can be managed only by the server, which in gradient approaches collects functions' gradients, averages them, does subsequent calculation to find a new point and broadcasts it to other devices. Such setting is considered to be popular, especially in theory. Another way of organizing communications is a decentralized setup, where all devices are connected in a single network with specific topology and, unlike centralized, can communicate with each other through networks' edges. In our work, we prefer the latter one.

The situation in which each worker's compressor is completely independent of the compressor operators of other devices seems to be rather improper for the decentralized setup. There are several reasons for such difficulties arising, the first is that with such approach we are not able to exploit data similarity between devices to the fullest extent. Explanation is simple, the difference between the averaging of compressed gradients and their uncompressed versions is expressed by the $AB$-inequality described in \citep{szlendak2021permutation}. For uncorrelated compressors, the value of constant $B$ is expressed as zero, while for correlated compressors it is significantly higher, which allows writing tighter estimates. The second problem we face lies on the level of technical implementation of communications in a network, the size of the transmitted packets in which is unknown in advance. Thus, it becomes complicated to perform such important operations as AllReduce \citep{basiccomm}, which is widely used in the present work, using arbitrary compressors. This is why correlated compressors, like PermK \citep{szlendak2021permutation}, are used in communication networks. In our paper we research both correlated and uncorrelated compressors, therefore to make algorithm suitable for different ways of organizing communications we change the name of specific communication operation with the word “broadcast”.  

\vspace{-0.2cm}
\subsection{Introduction of \texttt{DHPL-Katyusha}}
\vspace{-0.2cm}
At the beginning of the $k$-th iteration of \texttt{DHPL-Katyusha}, every worker computes $x^k$ locally, as a convex combination of $z^k, w^k$ and $y^k$ and broadcasts it to all other workers (line \ref{DHPL line1}). After this, each worker already has the same old and new point and therefore can compute the difference between local gradients in the current point $x^k$ and in the old point $w^k$ in the same manner (line \ref{DHPL line2}). Then, each worker applies a local compressor and sends this difference to everyone else (line \ref{DHPL line3}, line \ref{DHPL AllReduce}). Having done that, every device is able to find the average of all compressed differences and add full gradient in $w^k$ (stored locally on each device) to apply the variance reduction technique (line \ref{DHPL line4}). After that, “negative momentum” idea is used from \citep{kovalev2019dont} (line \ref{DHPL neg1}, line \ref{DHPL neg2}). Finally, $w^k$ is updated to the actual point, based on a coin flip, which all workers do with the same random seed (line \ref{DHPL line5}). If $w^k$ is adjusted, the uncompressed gradient is calculated by AllReduce procedure (line \ref{DHPL line6}). As $p$ strives to zero with $\omega \to \infty$, full gradient is updated rarely.

\begin{algorithm}[ht]
            \caption{\texttt{DHPL-Katyusha}}\label{alg:one}
            \textbf{Input:} initial $y^0 = w^0 = z^0 \in \mathbb{R}^d$, step size $\eta = \min\left\{\frac{\theta_2}{(1 + \theta_2) \theta_1}, \frac{\frac{\widetilde{L}}{L}\theta_2}{(1 + \theta_2) \theta_1} \right\}$, where $\widetilde{L} = L \frac{\omega}{n}$, $\sigma = \frac{\mu}{\widetilde{L}}$, parameters $\theta_1, \theta_2 \in \mathbb{R}$ and probability $p \in (0, 1]$ (every worker has the same random seed for calculating $p$).
            \begin{algorithmic}[1]
            \For {$k = 0, 1, 2,\dots K$} 
                \For {$i = 1 \dots n$ in parallel}
                    \vspace{0.4em}
                    \State $x^k \gets \theta_1 z^k + \theta_2 w^k + (1 - \theta_1 - \theta_2) y^k$ \label{DHPL line1}
                    \vspace{0.4em}
                    \State $g^k_i \gets \nabla f_{i} (x^k) - \nabla f_{i} (w^k)$
                    \vspace{0.5em}
                    \label{DHPL line2}
                    \State  $\widetilde{g}_i^k \gets Q_i (g_i^k)$ \label{DHPL line3}
                    \vspace{0.5em}
                    \State  Using communications broadcast  $\widetilde{g}_i^k$ \label{DHPL AllReduce}
                        \vspace{0.5em}
                    \State Compute $\widetilde{g}^k \gets \frac{1}{n} \sum\limits_{i = 1}^n \widetilde{g}_i^k + \nabla f(w^k)$\label{DHPL line4}
                    \vspace{0.5em}
                    \State $z^{k + 1} \gets \frac{1}{1 + \eta \sigma} (\eta \sigma x^k + z^k - \frac{\eta}{\widetilde{L}} \widetilde{g}^k )$\label{DHPL neg1}
                    \vspace{0.4em}
                    \State $y^{k + 1} \gets x^k + \theta_1 (z^{k + 1} - z^k)$ \label{DHPL neg2}
                    \State $w^{k + 1} \gets \displaystyle\begin{cases}
                         y^k, &\text{ with probability } p
                         \\
                         \label{DHPL line5}
                        w^k, &\text{ with probability }1-p
                    \end{cases}$
                    \vspace{0.5em}
                    \If{$w^{k + 1} \  \textbf{=} \ y^k$}
                        \vspace{0.4em}
                        \State Using communications broadcast  $\nabla f_i(w^{k + 1})$
                        \State Compute $\nabla f(w^{k + 1}) = \frac{1}{n}\sum\limits_{i = 1}^n \nabla f_i(w^{k + 1})$
                        \label{DHPL line6}
                    \EndIf
                \EndFor
            \EndFor
            \end{algorithmic}
\end{algorithm}

Note that the $Q_i$ compressors in this algorithm work independently of each other, only if the PermK operator is not considered.

\vspace{-0.2cm}
\subsection{Convergence results}
\vspace{-0.2cm}
To prove the asymptotics of \texttt{DHPL-Katyusha} we need to make the following assumptions for the problem \eqref{eq:hor}.
\begin{assumption} \label{ass:ass001}
Functions $f_{i}:\mathbb{R}^d \to \mathbb{R}$ are $L$-smooth for some $L > 0, \forall i \in \overline{1, n}$: 

$f_{i} (y) \leq f_{i} (x) + \la\nabla f_{i} (x), y - x\ra + \frac{L}{2} \|y - x\|^2 \text{ , $\forall x, y \in \mathbb{R}^d$}$.
\end{assumption}

Note that with such assumption, the target function $f(x)$ is also $L$ - smooth.

\begin{assumption} \label{ass:ass002}
The function $f :\mathbb{R}^d \to \mathbb{R}$ is $\mu$ - strongly convex for some $\mu > 0$:

$f (y) \geq f (x) + \la\nabla f (x), y - x\ra + \frac{\mu}{2} \| y - x \|^2 \text{ , $\forall x, y \in \mathbb{R}^d$}$.
\end{assumption}
In the original paper \citep{kovalev2019dont}, one of the steps in the proof of the asymptotics of \texttt{L-Katyusha}, which is called as Lemma 6, is the estimation of a term $\|g^k - \nabla f(x^k)\|^2$. It can be formally written in the following form:

\begin{lemma}
For the original \texttt{L-Katyusha} algorithm, the norm of the difference between the true and the real gradient can be estimated as
\begin{align*}
\|g^k - \nabla f(x^k)\|^2 \leq  2 L \left ( f (w^k) - f(x^k) - \langle \nabla f (x^k), w^k - x^k \rangle \right ).
\end{align*}
\end{lemma}

However, in our paper, $g^k$ is changed to $\widetilde{g}^k$, and therefore, it takes to utilize additional inequalities to gain similar terms for an upper estimate of $\|\widetilde{g}^k - \nabla f(x^k)\|^2$.

\begin{lemma}
For the \texttt{DHPL-Katyusha} algorithm, the norm of the difference between the true and the real gradient can be estimated as \label{lemm:lemma2.4}
\begin{align*}
\|\widetilde{g}^k - \nabla f(x^k)\|^2 \leq 2 \widetilde{L} \left ( f (w^k) - f(x^k) - \la \nabla f (x^k), w^k - x^k \ra \right ),
\end{align*}
where
\begin{align*}
\widetilde{L} = L\frac{\omega}{n}.
\end{align*}
\end{lemma}

As similar inequality is proven for several algorithms, we give the following definition:
\begin{definition}
In Katyusha-based algorithm we call $\widetilde{L}$ the efficient Lipschitz constant, if it is the smallest constant, such that:
\begin{align*}
\|\widetilde{g}^k - \nabla f(x^k)\|^2 \leq 2 \widetilde{L} \left ( f (w^k) - f(x^k) - \langle \nabla f (x^k), w^k - x^k \rangle \right ).
\end{align*}
\end{definition}
Note, that for the case $\widetilde{L} \geq L$ we can use all Lipschitz inequalities that were used in \citep{kovalev2019dont}, but for $\widetilde{L}$ and therefore conclude a proof. But not always this property is satisfied, and we need to use another inequalities to conclude a proof. Using this fact, we can formulate the convergence theorem:
\begin{theorem}
\label{th:one}
Let Assumptions~\ref{ass:ass001},~\ref{ass:ass002} be hold. Denote $\widetilde{L}$ as $L\frac{\omega}{n}$ and $x^*$ as the solution for the problem \ref{eq:hor}. Then after $k$ iterations of \texttt{DHPL-Katyusha}
\begin{align*}
\mathbb{E}\left[\mathbb{Z}^{k + 1} + \mathbb{Y}^{k + 1} + \mathbb{W}^{k + 1}\right]& \\
\leq \frac{1}{1 + \eta \sigma}\mathbb{Z}^k + (1 - \theta_1(1 -& \theta_2))\mathbb{Y}^k + \left(1 - \frac{p\theta_1}{1 + \theta_1} \right)\mathbb{W}^k,
\end{align*}
where:
\begin{align*}
&\mathbb{Z}^k \eqdef \frac{\widetilde{L} (1 + \eta \sigma)}{2 \eta} \|z^k - x^*\|^2,
\\
&\mathbb{Y}^k \eqdef \frac{1}{\theta_1} \left (f(y^k) - f(x^*) \right),
\\
&\mathbb{W}^k \eqdef \frac{\theta_2 (1 + \theta_1)}{p \theta_1} \left( f(w^k) - f(x^*) \right).
\end{align*}
\end{theorem}

We choose the concrete $p$ with the aim of reducing the average amount of information sent in each communication procedure. At each iteration, on average, workers communicate at a cost of $\mathcal{O}\left(\frac{1}{\beta} + p\right)$. Therefore, $p$ should be equal to $\frac{1}{\beta}$ to get a gain in asymptotics. Other parameters we choose are similar to that of \texttt{L-Katyusha} \citep{kovalev2019dont}.
\begin{corollary}
Denote $\beta$ as a density coefficient of compressor operator $Q$ from Definition \ref{ass:ass000}. Let $p = \frac{1}{\beta}$, $\theta_1 = \min{\{ \sqrt{\frac{2\sigma \beta }{3}} , \frac{1}{2} \}}$, $\theta_2 = \frac{1}{2}$. Then after 
$K = \mathcal{O}\left(\left(\sqrt{\frac{L}{\mu}}\left(\sqrt{\beta\frac{\omega}{n}} + 1\right) + \beta\right) \log{\frac{1}{\varepsilon}}\right)$ iterations $\mathbb{E}\left[\Psi^K\right] \leq \varepsilon \Psi^0$, where, the Lyapunov function $\Psi^K$ is defined as $\Psi^K \eqdef \mathbb{Z}^K + \mathbb{Y}^K + \mathbb{W}^K$.\\
Total information, sent by \texttt{DHPL-Katyusha} with such parameters, is  $$\mathcal{O}\left(\left(\sqrt{\frac{L}{\mu}}\left(\sqrt{\frac{1}{\beta}\frac{\omega}{n}} + 1 \right) + 1\right) \log{\frac{1}{\varepsilon}}\right).$$
\end{corollary}

\subsection{RandK and PermK comparison}
For RandK compressor, the exact $\omega$ is equal to $\frac{d}{K}$ \citep{beznosikov2024biased}. Therefore, using Lemma \ref{lemm:lemma2.4} we get the efficient Lipschitz constant, which is always greater than $L$:

\begin{lemma}\label{theorem:RandK horizontal}
\textup{(The efficient Lipschitz for RandK compressor in the horizontal case)}\\
For \texttt{DHPL-Katyusha} with RandK compressor, the efficient  Lipschitz constant is less than:
\begin{equation*}
\widetilde{L}_{Rand} = \left(\frac{d}{nK} + 1 \right) L
\end{equation*}
\end{lemma}

As $Q_i$ for PermK operator are correlated, the analysis for it is slightly different from the original \texttt{DHPL-Katyusha}. But we still managed to find its efficient Lipschitz constant and therefore can estimate asymptotics for it:

\begin{lemma}\label{lem:PermK horizontal}
\textup{(The efficient Lipschitz for PermK compressor in the horizontal case)}\\
For \texttt{DHPL-Katyusha} with PermK compressor, the efficient Lipschitz constant is equal to $L$.
\end{lemma}

As it can be seen in the above lemmas, $\widetilde{L}_{Rand}$ and $\widetilde{L}_{Perm}$ are asymptotically equal as we take $K = \frac{d}{n}$ for RandK in order to compare with PermK. However, as it is shown in the experiments, PermK compressor is empirically better for horizontal data division case than RandK.

\vspace{-0.2cm}
\section{Vertical Case}\label{sec:zoo}
\vspace{-0.2cm}

In this section, we introduce the \texttt{Katyusha}-based algorithm of distributed optimization with vertical data division, which is named as \texttt{DVPL-Katyusha}. We introduce an algorithm for the problem of minimization of the linear loss function. This can be formally written as:

\begin{equation}
\label{eq:Arbitrary loss}
\min_{x \in \R^d} \left[\mathcal{L}(Ax, b) \eqdef \frac{1}{s} \sum\limits_{j = 1}^s \mathcal{L}_j(A_j^Tx, b_j) \right],
\end{equation}

where $b \in \mathbb{R}^s$ is a vector of targets, $A \in \R^{s \times d}$ is a feature matrix, $s$ is a number of samples and $d$ is a number of features. We denote $A_j$ as a $j$-th row of matrix $A$.

In the vertical setup, we assume that every worker has its own set of features, which is formally represented, as a set of columns of the feature matrix $A$. Thus, we can rewrite the dot product $A_j^Tx$ in the form of:
\begin{equation}
A_j^Tx = \sum\limits_{i = 1}^n A_{ji}^T x_i,
\end{equation}
where as $A_i$ and $x_i$ we denote submatrix and point subvector corresponding to the set of $i$-th worker's features, $n$ is a number of devices.

Compressor in \texttt{DVPL-Katyusha} is basically a RandK compressor \cite{beznosikov2024biased}, which originally compresses data by choosing random coordinates with certain probabilities and broadcasts information only from them. The main difference in our case if that every worker has the same random seed for such choice, therefore it can be viewed as an analogue to batching technique for the vertical case.

\vspace{-0.2cm}
\subsection{Introduction of \texttt{DVPL-Katyusha}}
\vspace{-0.2cm}
In \texttt{DVPL-Katyusha} every worker has its own set of features. Therefore, unlike the horizontal regime, in the vertical case all operations are performed on subvectors, corresponding to individual worker's components. We denote such subvectors with the same letter as the vector itself, but with lower index $i$. Similar to \texttt{DHPL-Katyusha} (Algorithm \ref{alg:one}) at each iteration every worker sets $x^k_i$ as a convex combination of $z^k_i, w^k_i$ and $y^k_i$ (line \ref{DVPL line1}). A gradient calculated on the $j$-th sample depends on the dot product of $A^T_{j}$ and a point in which we want to perform a computation. This is why, to perform a step, we need to obtain dot products $\langle A^T_j, x^k \rangle$ and $\langle A^T_j, w^k \rangle$ for every sample, chosen by RandK. We want to select samples with large Lipschitz constant more often, therefore, we select each sample with probability $\frac{L_j}{s\bar{L}}$.  As every worker can compute only its own part of this sum (in our denotations for a single sample $j$ it is $\langle A^T_{ji}, x^k_{i} \rangle$ or $\langle A^T_{ji}, w^k_{i} \rangle$), all such parts from every chosen sample should be broadcasted in order to compute $g^k$ (line \ref{DVPL compute}, line \ref{DVPL AllReduce}). After this, a standard momentum part of \texttt{Katyusha} is performed locally (line \ref{DVPL line3}, line \ref{DVPL negative}), after which $w^k$ is updated with some probability, with every worker having the same random seed for it (line \ref{DVPL line4}). If $w^k$ is adjusted to the new point, the full gradient (e.g., all samples are used) in this vector is calculated the same way as is written above, but without any compression (line \ref{DVPL line5}). Note that similar to the horizontal case, we consider probability of such update low, e.g., full gradient is calculated infrequently. 

\begin{algorithm}
            \caption{\texttt{DVPL-Katyusha}} \label{alg:orig}
            \textbf{Input:} initial $y^0 = w^0 = z^0 \in \mathbb{R}^d$, step size $\eta = \frac{\theta_2}{(1 + \theta_2) \theta_1}$, $\sigma = \frac{K\mu}{\bar{L}}$ if $\frac{\bar{L}}{K} \geq L$ or else $\sigma = \frac{\mu}{L}$, parameters $\theta_1, \theta_2 \in \mathbb{R}$ and probability $p \in (0, 1]$, every worker has the same random seed for RandK random. RandK select $j$-th sample with probability $p_j = \frac{L_j}{s\bar{L}}$.
            \vspace{0.3em}
            \begin{algorithmic}[1]
            \For {$k = 0, 1, 2,\dots K$}
            \vspace{0.4em}
            \For {$i = 1 \dots n$ in parallel}
                \vspace{0.2em}
                \State $x^k_{i} \gets \theta_1 z^k_{i} + \theta_2 w^k_{i} + (1 - \theta_1 - \theta_2) y^k_{i}$ \label{DVPL line1}
                \vspace{0.2em}
                \State Compute $X_i^k=\text{RandK}\left(\left\|\left<A^T_{ji}, x^k_{i}\right> \right \|_{j=\overline{1,s}} \right)$ 
                \vspace{0.2em} \label{DVPL compute}
                \State Compute $W_i^k=\text{RandK}\left(\left\|\left<A^T_{ji}, w^k_{i}\right> \right \|_{j=\overline{1,s}} \right)$
                \vspace{0.2em}
                \State  Using communications broadcast $X_i^k$ and $W_i^k$ \label{DVPL AllReduce}
                \vspace{0.4em}
                \State $J^k = \{j^k_1\!, \cdots\!, j^k_n \}$ - indices, selected by RandK
                \vspace{0.3em}
                \State $g^k_{i}\gets\frac{1}{K}\sum\limits_{j \in {J^k}}\frac{1}{sp_j}\nabla \mathcal{L}_j\left(\sum\limits_{i = 1}^{n}X^k_{ij}, b_j\right)_i - \frac{1}{K}\sum\limits_{j \in {J^k}}\frac{1}{sp_j}\nabla \mathcal{L}_j\left(\sum\limits_{i = 1}^{n}W^k_{ij}, b_j\right)_i +  \nabla\mathcal{L}\left(Aw^k, b\right)_i$ 
                \vspace{0.2cm}
                \State $z^{k + 1}_{i} \gets \frac{1}{1 + \eta \sigma} (\eta \sigma x^k_{i} + z^k_{i} - \frac{\eta}{\widetilde{L}} g^k_{i} )$  \label{DVPL line3}
                \vspace{0.2cm}
                \State $y^{k + 1}_{i} \gets x^k_{i} + \theta_1 (z^{k + 1}_{i} - z^k_{i})$  \label{DVPL negative}
                \vspace{0.2cm}
                \State $w^{k + 1}_{i} \gets \displaystyle\begin{cases} \label{DVPL line4}
                     y^k_{i}, &\text{ with probability } p\\ \vspace{-0.9em}\\
                    w^k_{i}, &\text{ with probability }1-p
                \end{cases}
                $
                \vspace{0.5em}
                \If{$w^{k + 1}_{i} \ \textbf{=} \  y^k_{i}$}
                    \For {$j = 1\dots s$}
                        \State Compute $ \left<A^T_{ji}, w^k_{i}\right>$ \label{DVPL line5}
                        \State Using communications broadcast $\left<A^T_{ji}, w^k_{i}\right>$
                    \EndFor
                    \State Compute $\nabla\mathcal{L}\left(Aw^k, b\right)_i$
                \EndIf
            \EndFor 
        \EndFor
            \end{algorithmic}
\end{algorithm}

\subsection{Convergence results}

As in \texttt{DHPL-Katyusha} we need to assume our target function to be strong convex and worker's functions to be $L$-smooth.
\begin{assumption} \label{ass: smoothness DVPL}
Functions $\mathcal{L}_{i}:\mathbb{R}^d \to \mathbb{R}$ are $L$-smooth for some $L > 0, \forall i \in \overline{1, s}$: 

$\mathcal{L}_{i} (y)\!\leq\!\mathcal{L}_{i} (x)\!+\!\la\nabla \mathcal{L}_{i} (x), y\!-\!x\ra + \frac{L}{2} \|y - x\|^2 \text{ , $\forall x, y \in \mathbb{R}^d$}$.
\end{assumption}

\begin{assumption} \label{ass: strong convex for DVPL}
The function $\mathcal{L} :\mathbb{R}^d \to \mathbb{R}$ is $\mu$ - strongly convex for some $\mu > 0$:

$\mathcal{L} (y) \geq \mathcal{L} (x) + \la\nabla \mathcal{L} (x), y - x\ra + \frac{\mu}{2} \| y - x \|^2 \text{ , $\forall x, y \in \mathbb{R}^d$}$.
\end{assumption}

To estimate the asymptotics of the algorithm, we need to get a bound of the efficient Lipschitz constant for \texttt{DVPL-Katyusha}. Denote $\bar{L}$ as $\frac{1}{s} \sum \limits_{i = 1}^s L_i$. Therefore, we can formulate the following lemma.
\begin{lemma} \label{lem: lemma 6 arbitrary loss}
The efficient Lipschitz constant for \texttt{DVPL-Katyusha} is less than $\widetilde{L}$, where:
\begin{align*}
\widetilde{L} = \max \left\{L, \frac{\bar{L}}{K} \right\}.
\end{align*}
\end{lemma}

In this case $L$ if always less than $\widetilde{L}$, therefore knowing the efficient Lipschitz constant, we can estimate a convergence rate of our algorithm:

\begin{theorem}
\label{th:bad_theorem}
Let assumptions \ref{ass: smoothness DVPL} and \ref{ass: strong convex for DVPL} be hold. Denote $\widetilde{L}$ as $\max \left\{L, \frac{\bar{L}}{K} \right\}$ and $x^*$ as the solution for the problem \ref{eq:Arbitrary loss}. Then after $k$ iterations of \texttt{DVPL-Katyusha}
\begin{align*}
\mathbb{E}\left[\mathbb{Z}^{k + 1} + \mathbb{Y}^{k + 1} + \mathbb{W}^{k + 1}\right] 
\leq \frac{1}{1 + \eta \sigma}\mathbb{Z}^k + (1 - \theta_1(1 -& \theta_2))\mathbb{Y}^k + \left(1 - \frac{p\theta_1}{1 + \theta_1} \right)\mathbb{W}^k,
\end{align*}\\
where:
\begin{align*}
&\mathbb{Z}^k \eqdef \frac{\widetilde{L} (1 + \eta \sigma)}{2 \eta} \|z^k - x^*\|^2,
\\
&\mathbb{Y}^k \eqdef \frac{1}{\theta_1} \left (f(y^k) - f(x^*) \right),
\\
&\mathbb{W}^k \eqdef \frac{\theta_2 (1 + \theta_1)}{p \theta_1} \left( f(w^k) - f(x^*) \right).
\end{align*}
\end{theorem}

After choosing concrete parameters, we can get an estimation for a total number of iterations and sent information.

\begin{corollary}
Let $p = \frac{K}{s}$, $\theta_1 = \min{\{ \sqrt{\frac{2\sigma sK }{3}} , \frac{1}{2} \}}$, $\theta_2 = \frac{1}{2}$. If $\frac{\bar{L}}{K} \geq L$, then after $N = \mathcal{O}\left(\left(\sqrt{\frac{s\bar{L}}{\mu}} + s \right) \log{\frac{1}{\varepsilon}}\right)$ iterations $\mathbb{E}\left[\Psi^N\right] \leq \varepsilon \Psi^0$, where the Lyapunov function $\Psi^N$ is defined as $\Psi^N \eqdef \mathbb{Z}^N + \mathbb{Y}^N + \mathbb{W}^N$. If $L > \frac{\bar{L}}{K}$, then after $\mathcal{O}\left(\left(\sqrt{\frac{sKL}{\mu}} + s \right) \log{\frac{1}{\varepsilon}}\right)$ iterations, the same accuracy is achieved.
\end{corollary}

The average information cost of a single iteration is $\mathcal{O}\left(K + ps\right)$. Therefore, we choose $p$ as $\frac{K}{s}$ to asymptotically reduce it.
\subsection{Scalar compressors}
This section introduces \texttt{DVPL-Katyusha with scalar compression}, which is the generalization of Algorithm \ref{alg:orig} with the compressors that operate on the workers' parts of the dot product $A^T_{ji}x^k_i$. For simplicity, only the MSE loss function will be considered in this section, which can be formally expressed as:
\begin{equation}
\label{eq:MSE}
\min_{x \in \R^d} \left[f(x) \eqdef \frac{1}{s} ||Ax - b||^2 \right].
\end{equation}
As our setup is still vertical, we can rewrite the above in terms of workers' components:
\begin{equation}
\min_{x \in \R^d} \left[f(x) \eqdef \frac{1}{s}\left\| \sum\limits_{i = 1}^n A_ix_i - b \right \|^2 \right].
\end{equation}

The practical utility of the MSE loss function is significant, as linear models are quick to train and can serve as effective benchmarks. Additionally, large models such as neural networks can leverage a trained model to generate features and then utilize only the final linear layer. Furthermore, implementing privacy mechanisms in vertical federated learning is simpler with this approach through linear models, like in \citep{huang2022coresets}.

\begin{algorithm}
            \caption{\texttt{DVPL-Katyusha with scalar compression}}\label{alg:two}
            \textbf{Input:} initial $y^0 = w^0 = z^0 \in \mathbb{R}^d$,  step size $\eta = \frac{\theta_2}{(1 + \theta_2) \theta_1}$, $\sigma = \frac{\mu}{\widetilde{L}}$, where $\textstyle{\widetilde{L} = L\left(1 + (\omega - 1)\frac{s \cdot \sum\limits_{j = 1}^s L_j^2}{\mu^2} \right)}$, parameters $\theta_1, \theta_2 \in \mathbb{R}$ and probability $p \in (0, 1]$, every worker has the same random seed for RandK random. RandK select $j$-th sample with probability $p_j = \frac{1}{s}$.
            \vspace{0.2em}
            \begin{algorithmic}[1]
            \For {$k = 0, 1, 2,\dots K$}
            \vspace{0.4em}
            \For {$i = 1 \dots n$ in parallel}
                \vspace{0.2em}
                \State $x^k_{i} \gets \theta_1 z^k_{i} + \theta_2 w^k_{i} + (1 - \theta_1 - \theta_2) y^k_{i}$ 
                \vspace{0.2em}
                \State Compute $D^k_i=\text{RandK}\left(\left\|\left<A^T_{ji},\!x^k_{i}\!-\!w^k_i\right> \right \|_{j=\overline{1,s}} \right)$ 
                \vspace{0.2em}
                \State $J^k = \{j^k_1\!, \cdots\!, j^k_n \}$ - indices, selected by RandK
                \vspace{0.4em}
                \State  Using communications broadcast $Q_i\left(D^k_i\right)$
                \vspace{0.3em}
                \State $g^k_{i}\gets\frac{2}{b_s}\sum\limits_{j \in {J^k}}A^T_{ji}\sum\limits_{i = 1}^{n}Q_i\left(D_{ij}\right) + \frac{2}{s}\left(A^TAw^k-A^Tb\right)_{i}$
                \vspace{0.2cm}
                \State $z^{k + 1}_{i} \gets \frac{1}{1 + \eta \sigma} (\eta \sigma x^k_{i} + z^k_{i} - \frac{\eta}{\widetilde{L}} g^k_{i} )$ 
                \vspace{0.2cm}
                \State $y^{k + 1}_{i} \gets x^k_{i} + \theta_1 (z^{k + 1}_{i} - z^k_{i})$ 
                \vspace{0.2cm}
                \State $w^{k + 1}_{i} \gets \displaystyle\begin{cases} 
                     y^k_{i}, &\text{ with probability } p\\ \vspace{-0.9em}\\
                    w^k_{i}, &\text{ with probability }1-p
                \end{cases}
                $
                \vspace{0.5em}
                \If{$w^{k + 1}_{i} \ \textbf{=} \  y^k_{i}$}
                    \For {$j = 1\dots s$}
                        \State Compute $ \left<A^T_{j_{i}}, x^k_{i} - w^k_{i} \right>$
                        \State Using communications  broadcast $ \left<A^T_{ji}, x^k_{i} - w^k_{i} \right>$
                    \EndFor
                    \State Compute $\frac{2}{s}\left(A^TAw^k\!-\!A^Tb\right)_{i}$ 
                \EndIf
            \EndFor 
        \EndFor
            \end{algorithmic}
\end{algorithm}

Like in Algorithm \ref{alg:orig} we need to assume $L$-smoothness and $\mu$ strong convexity of our problem. In the MSE loss, this can be written in terms of eigenvalues.
\begin{assumption} \label{ass:ass005}
The matrix $A^TA$ has all eigenvalues bounded in a segment of $\left[\mu, L \right]$.
\end{assumption}

To estimate the asymptotics of the algorithm, we need to get a bound of the efficient Lipschitz constant for \texttt{DVPL-Katyusha with scalar compression}. 

\begin{lemma} \label{lem: lemma 6 vertically distributed} 
The efficient Lipschitz constant for \texttt{DVPL-Katyusha with scalar compression} is less than $\widetilde{L}$, where:
\begin{align*}
\widetilde{L} = L\left(1+(\omega-1)\frac{s\cdot \sum\limits_{j = 1}^s L_j^2}{\mu^2} \right).
\end{align*}
\end{lemma}
We also proved that under made assumptions there exists a function that has $\widetilde{L}$ proportional to $\frac{L^3}{\mu^2}$, therefore our asymptotics for \texttt{DVPL-Katyusha with scalar compression} is optimal in the terms of $L$ and $\mu$ constants.

\begin{lemma} \label{lem: lower bound vertical}
There exists such function that holds under \ref{ass:ass005} that its efficient Lipschitz constant $\widetilde{L}$ is proportional to  $(\omega - 1) \frac{L^3}{\mu^2}$.
\end{lemma}

\subsection{PermK compressor}
Beside RandK and scalar operators, we have also implemented another compression strategy \citep{szlendak2021permutation}. In PermK every worker choose its own set of samples, which can not intersect with any other set. After that, each worker do AllReduce with calculated values. Therefore, we can prove the following lemma:
\newpage
\begin{lemma} \label{lem: vertical PermK}
The efficient Lipschitz constant for PermK compressor is less than:
\begin{align*}
\widetilde{L}_{Perm} = {2 L} \frac{s\sum\limits_{j = 1}^s L_j^2}{\mu^2}.
\end{align*}
\end{lemma}

Unlike horizontal regime, RandK has a superior theoretical asymptotics on PermK.

\begin{figure*}[t!]

\begin{minipage}{0.33\textwidth}
  \centering
\includegraphics[width =  0.8\textwidth ]{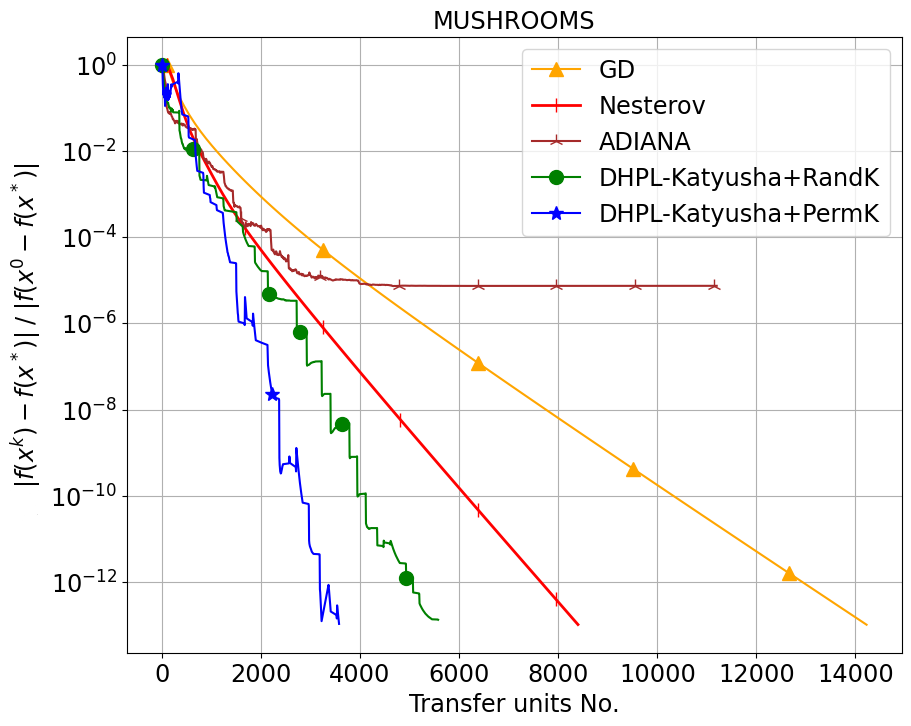}
\end{minipage}%
\begin{minipage}{0.33\textwidth}
  \centering
\includegraphics[width =  0.8\textwidth ]{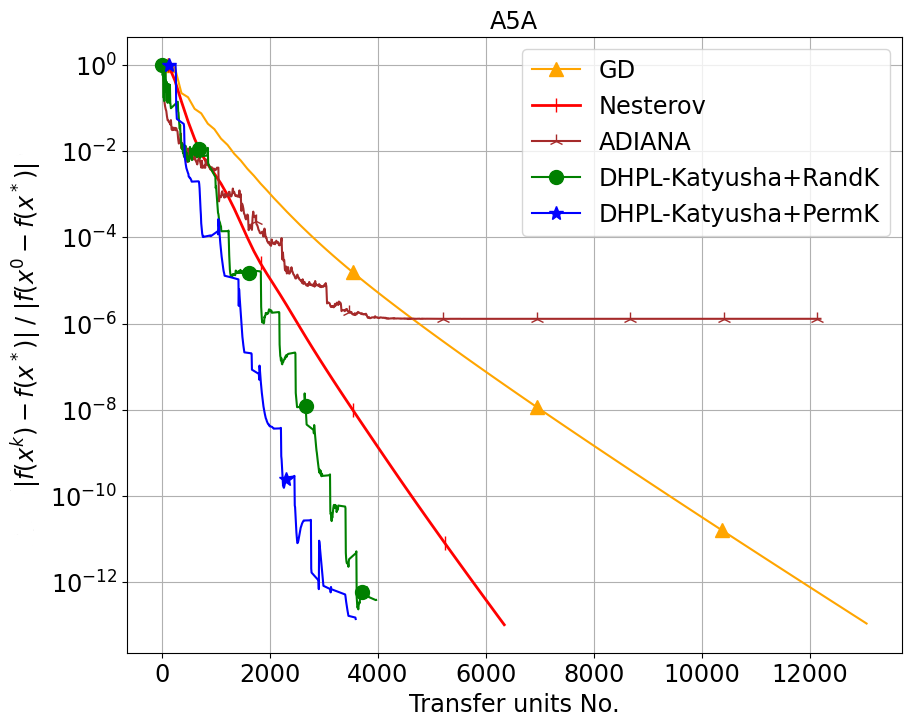}
\end{minipage}%
\begin{minipage}{0.33\textwidth}
  \centering
\includegraphics[width =  0.8\textwidth ]{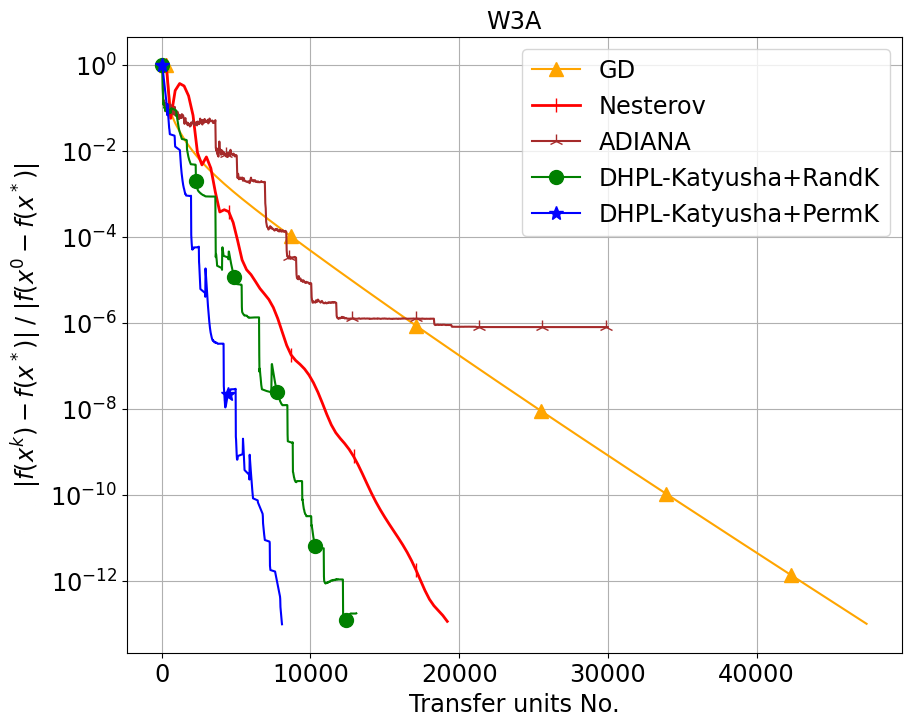}
\end{minipage}%
\\
\begin{minipage}{0.01\textwidth}
\quad
\end{minipage}%
\begin{minipage}{0.33\textwidth}
  \centering
(a) \texttt{mushrooms}
\end{minipage}%
\begin{minipage}{0.33\textwidth}
\centering
 (b) \texttt{a5a}
\end{minipage}%
\begin{minipage}{0.33\textwidth}
\centering
  (c) \texttt{w3a}
\end{minipage}%

\caption{
Comparison of different algorithms for solving optimization problem from section \ref{subsection4.1} in horizontal case on LIBSVM datasets \texttt{mushrooms}, \texttt{a5a} and \texttt{w3a}.}
    \label{fig:comparison1}
    
\end{figure*}
\begin{figure*}[t!]

\begin{minipage}{0.33\textwidth}
  \centering
\includegraphics[width =  0.8\textwidth ]{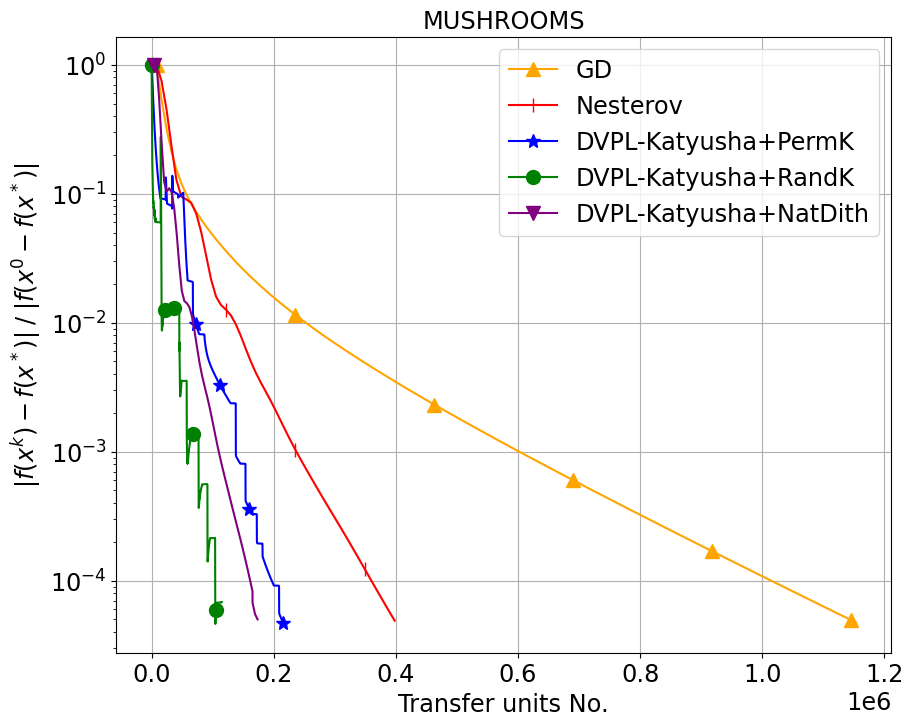}
\end{minipage}%
\begin{minipage}{0.33\textwidth}
  \centering
\includegraphics[width =  0.8\textwidth ]{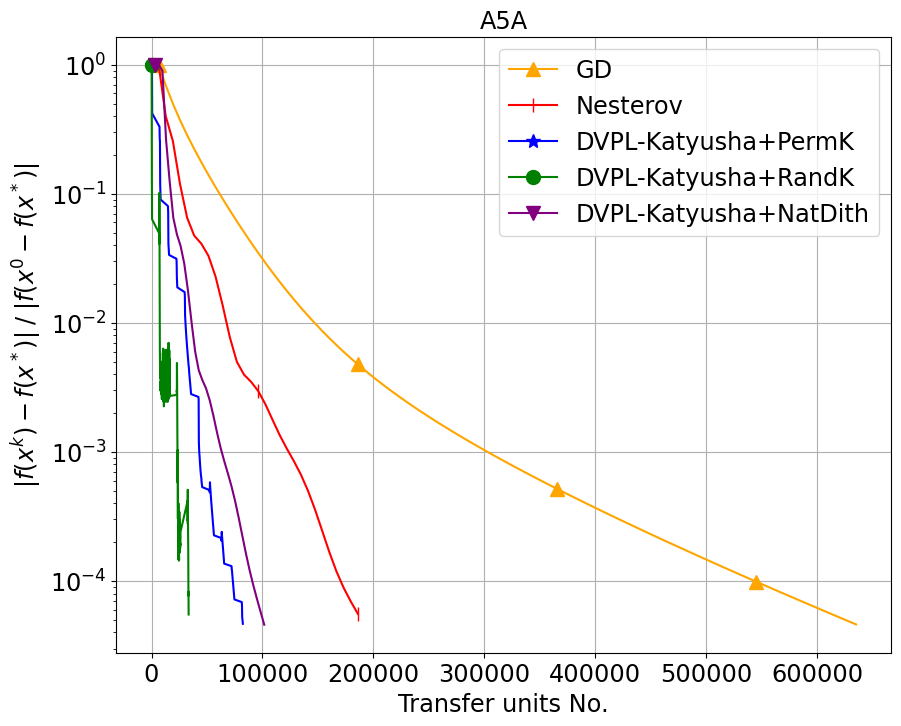}
\end{minipage}%
\begin{minipage}{0.33\textwidth}
  \centering
\includegraphics[width =  0.8\textwidth ]{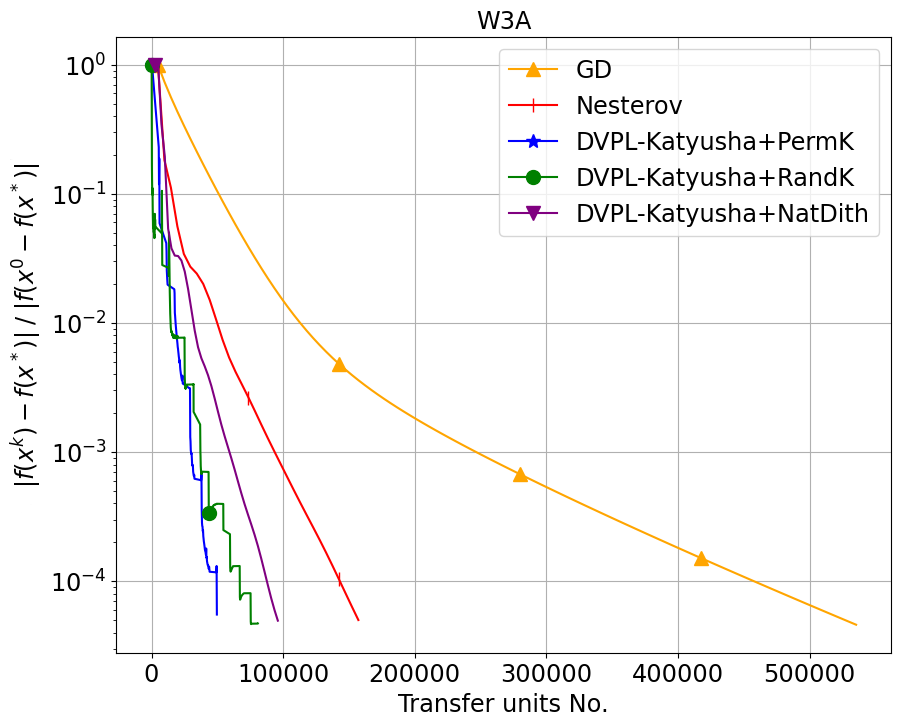}
\end{minipage}%
\\
\begin{minipage}{0.01\textwidth}
\quad
\end{minipage}%
\begin{minipage}{0.33\textwidth}
  \centering
(a) \texttt{mushrooms}
\end{minipage}%
\begin{minipage}{0.33\textwidth}
\centering
 (b) \texttt{a5a}
\end{minipage}%
\begin{minipage}{0.33\textwidth}
\centering
  (c) \texttt{w3a}
\end{minipage}%

\caption{
Comparison of different algorithms for solving optimization problem from section \ref{subsection4.2} in vertical case on LIBSVM datasets \texttt{mushrooms}, \texttt{a5a} and \texttt{w3a}.}
    \label{fig:comparison2}
    \vskip-12pt
\end{figure*}

\section{Numerical Experiments}
In this section, we conduct multiple numerical experiments on horizontal and vertical cases with a binary classification problem and a strongly convex target function on \texttt{mushrooms}, \texttt{a5a} and \texttt{w3a} datasets from the LIBSVM library.
\subsection{Binary classification in the horizontal regime} \label{subsection4.1}
As assumed in our paper, in horizontal case we present our target function as the following finite sum:
\begin{equation*}
\min_{x \in \mathbb{R}^d} \left\{ f(x) 
= \frac{1}{n}\sum\limits_{m = 1}^n f_m(x)\right\},
\end{equation*}
where $f_m(x) = \frac{1}{s / n} \sum_{i = 1}^{s / n}\log(1 + \exp(-y_i a_i^T x)) + \frac{\lambda}{2} \|x\|^2_2$ -- the sum here is taken with samples ($a_i$, $y_i$) from the local dataset of the $m$-th worker, $\lambda = L / 100$ -- $L$ here is the target function's $f(x)$ smoothness constant.

In our experiments on horizontal case, we have $n = 100$ devices involved in the training process and compare \texttt{DHPL-Katyusha} in combination with Rand1\% and PermK compressors with accelerated algorithm \texttt{ADIANA} and vanilla GD and AGD, running in the same distributed setting. All of the plots on the Figure \ref{fig:comparison1} are constructed for the convergence with tuned parameters. It is known from the article \cite{basiccomm} that the AllReduce of vectors of size $d$ takes time proportional to the magnitude of $2 \frac{n-1}{n} d$, where d is the size of the vectors transmitted in this way. On the $x$-axis, we postponed either $\frac{d}{\beta}$ or $d$, depending on the outcome of the coinflip produced during the iteration. Thus, on the x-axis, we have a number proportional to the number of numpy.doubles (without taking into account the $2 \frac{n-1}{n}$ multiplier that is the same for all the plots) transmitted during the entire learning process.

As it was expected in theory, our algorithm has an advantage over AGD with a total communication complexity of $\mathcal{O}\left(\sqrt{\frac{L}{\mu}}\right)$ due to the $n$ term in the denominator of our algorithm's best complexity -- $\mathcal{O}\left(\sqrt{\frac{L}{n \mu}}\right)$. More than that, PermK clearly performs better than Rand1\% in such setup, as can be seen from the graph.

Finally, it worth to be noted that the ADIANA algorithm in our setup could not be run to convergence, the algorithm converges only to the neighborhood of true solution, which is not stated in theory. At the same time, it is worth noting that a wide range of algorithm parameters were considered and some of its modifications were considered too in order to run it to the full convergence.
\subsection{Linear regression in the vertical regime} \label{subsection4.2}

In the case of vertical data partition, again, we present our target function as:
\begin{equation*}
\min_{x \in \mathbb{R}^d} \left\{ f(x) 
= \frac{1}{s}\sum\limits_{i = 1}^s l_i(a_i^Tx, y_i) + \frac{\lambda}{2} \|x\|^2_2\right\},
\end{equation*}
where $a_i \in \mathbb{R}^d, y_i \in \{ -1, 1 \}, \forall i \in 1, \dots, s$, $\lambda = L / 100$. $l_i : \mathbb{R}^2 \to \mathbb{R}$: $l_i(a_i^T x, y_i) = (a_i^T x - y_i)^2$.

We also conducted experiments to compare \texttt{DVPL-Katyusha} in combination with Rand1\%, Natural Dithering \cite{horváth2019stochastic} and PermK with vertical versions of regular gradient descent and its Nesterov acceleration in the same distributed setting. Number of workers in the experiments we conducted is $n = 5$. And again, on the x-axis, we postpone either $\frac{b}{\beta}$ or $s$ for each iteration, without taking into account the coefficient $2\frac{n-1}{n}$, which is the same for all implementations.

Overall, it is clear from the graph that there is a significant acceleration of convergence compared to the naive implementation of the Nesterov approach in a vertically distributed mode (when devices share all $s$ values of their local models, after which they calculate their part of the gradient and produce a Nesterov descent based on it) and with the same implementation of ordinary gradient descent.

In addition, the superiority of the PermK compressor is no longer present in this case, and the Rand1\% compressor comes out on top, which emphasizes the different importance of compressors' properties for different data division modes.

\section*{Impact Statement}
This paper presents work whose goal is to advance the field of Machine Learning. There are many potential societal consequences of our work, none of which we feel must be specifically highlighted here.

%%%%%%%%%%%%%%%%%%%%%%%%%%%%%%%%%%%%%%%%%%%%%%%%%%%%%%
\bibliography{ref}

\begin{thebibliography}{34}
\providecommand{\natexlab}[1]{#1}
\providecommand{\url}[1]{\texttt{#1}}
\expandafter\ifx\csname urlstyle\endcsname\relax
  \providecommand{\doi}[1]{doi: #1}\else
  \providecommand{\doi}{doi: \begingroup \urlstyle{rm}\Url}\fi

\bibitem[Alacaoglu and Malitsky(2022)]{alacaoglu2022stochastic}
A.~Alacaoglu and Y.~Malitsky.
\newblock Stochastic variance reduction for variational inequality methods,
  2022.

\bibitem[Alistarh et~al.(2017)Alistarh, Grubic, Li, Tomioka, and
  Vojnovic]{alistarh2017qsgd}
D.~Alistarh, D.~Grubic, J.~Li, R.~Tomioka, and M.~Vojnovic.
\newblock Qsgd: Communication-efficient sgd via gradient quantization and
  encoding, 2017.

\bibitem[Allen-Zhu(2018)]{BadKatyusha}
Z.~Allen-Zhu.
\newblock Katyusha: The first direct acceleration of stochastic gradient
  methods.
\newblock \emph{Journal of Machine Learning Research}, 18\penalty0
  (221):\penalty0 1--51, 2018.

\bibitem[Beznosikov and Gasnikov(2023)]{beznosikov2023similarity}
A.~Beznosikov and A.~Gasnikov.
\newblock Similarity, compression and local steps: Three pillars of efficient
  communications for distributed variational inequalities, 2023.

\bibitem[Beznosikov et~al.(2023)Beznosikov, Richtárik, Diskin, Ryabinin, and
  Gasnikov]{beznosikov2023distributed}
A.~Beznosikov, P.~Richtárik, M.~Diskin, M.~Ryabinin, and A.~Gasnikov.
\newblock Distributed methods with compressed communication for solving
  variational inequalities, with theoretical guarantees, 2023.

\bibitem[Beznosikov et~al.(2024)Beznosikov, Horváth, Richtárik, and
  Safaryan]{beznosikov2024biased}
A.~Beznosikov, S.~Horváth, P.~Richtárik, and M.~Safaryan.
\newblock On biased compression for distributed learning, 2024.

\bibitem[Cai et~al.(2022)Cai, Fan, Kang, Fan, Xu, Wang, and Yang]{Cai_2022}
D.~Cai, T.~Fan, Y.~Kang, L.~Fan, M.~Xu, S.~Wang, and Q.~Yang.
\newblock Accelerating vertical federated learning.
\newblock \emph{IEEE Transactions on Big Data}, page 1–10, 2022.
\newblock ISSN 2372-2096.
\newblock \doi{10.1109/tbdata.2022.3192898}.
\newblock URL \url{http://dx.doi.org/10.1109/TBDATA.2022.3192898}.

\bibitem[Castiglia et~al.(2023)Castiglia, Das, Wang, and
  Patterson]{castiglia2023compressedvfl}
T.~Castiglia, A.~Das, S.~Wang, and S.~Patterson.
\newblock Compressed-vfl: Communication-efficient learning with vertically
  partitioned data, 2023.

\bibitem[Chan et~al.(2007)Chan, Heimlich, Purkayastha, and Van
  De~Geijn]{basiccomm}
E.~Chan, M.~Heimlich, A.~Purkayastha, and R.~Van De~Geijn.
\newblock Collective communication: theory, practice, and experience.
\newblock \emph{Concurrency and Computation: Practice and Experience},
  19\penalty0 (13):\penalty0 1749--1783, 2007.

\bibitem[Chen et~al.(2021)Chen, Ma, Fan, Kang, Xu, and
  Yang]{chen2021secureboost}
W.~Chen, G.~Ma, T.~Fan, Y.~Kang, Q.~Xu, and Q.~Yang.
\newblock Secureboost+ : A high performance gradient boosting tree framework
  for large scale vertical federated learning, 2021.

\bibitem[Chilimbi et~al.(2014)Chilimbi, Suzue, Apacible, and
  Kalyanaraman]{chilimbi2014project}
T.~Chilimbi, Y.~Suzue, J.~Apacible, and K.~Kalyanaraman.
\newblock Project adam: Building an efficient and scalable deep learning
  training system.
\newblock In \emph{11th USENIX symposium on operating systems design and
  implementation (OSDI 14)}, pages 571--582, 2014.

\bibitem[Gorbunov(2021)]{gorbunov2021distributed}
E.~Gorbunov.
\newblock Distributed and stochastic optimization methods with gradient
  compression and local steps, 2021.

\bibitem[Gorbunov et~al.(2022)Gorbunov, Burlachenko, Li, and
  Richtárik]{marina}
E.~Gorbunov, K.~Burlachenko, Z.~Li, and P.~Richtárik.
\newblock Marina: Faster non-convex distributed learning with compression,
  2022.

\bibitem[Goyal et~al.(2018)Goyal, Dollár, Girshick, Noordhuis, Wesolowski,
  Kyrola, Tulloch, Jia, and He]{goyal2018accurate}
P.~Goyal, P.~Dollár, R.~Girshick, P.~Noordhuis, L.~Wesolowski, A.~Kyrola,
  A.~Tulloch, Y.~Jia, and K.~He.
\newblock Accurate, large minibatch sgd: Training imagenet in 1 hour, 2018.

\bibitem[Gu et~al.(2020)Gu, Xu, Huo, Deng, and Huang]{gu2020privacypreserving}
B.~Gu, A.~Xu, Z.~Huo, C.~Deng, and H.~Huang.
\newblock Privacy-preserving asynchronous federated learning algorithms for
  multi-party vertically collaborative learning, 2020.

\bibitem[He et~al.(2024)He, Huang, and Yuan]{he2024unbiased}
Y.~He, X.~Huang, and K.~Yuan.
\newblock Unbiased compression saves communication in distributed optimization:
  When and how much?
\newblock \emph{Advances in Neural Information Processing Systems}, 36, 2024.

\bibitem[Horváth et~al.(2019)Horváth, Kovalev, Mishchenko, Stich, and
  Richtárik]{horváth2019stochastic}
S.~Horváth, D.~Kovalev, K.~Mishchenko, S.~Stich, and P.~Richtárik.
\newblock Stochastic distributed learning with gradient quantization and
  variance reduction, 2019.

\bibitem[Huang et~al.(2022)Huang, Li, Sun, and Zhao]{huang2022coresets}
L.~Huang, Z.~Li, J.~Sun, and H.~Zhao.
\newblock Coresets for vertical federated learning: Regularized linear
  regression and $k$-means clustering, 2022.

\bibitem[Johnson and Zhang(2013)]{NIPS2013_ac1dd209}
R.~Johnson and T.~Zhang.
\newblock Accelerating stochastic gradient descent using predictive variance
  reduction.
\newblock In C.~Burges, L.~Bottou, M.~Welling, Z.~Ghahramani, and
  K.~Weinberger, editors, \emph{Advances in Neural Information Processing
  Systems}, volume~26. Curran Associates, Inc., 2013.
\newblock URL
  \url{https://proceedings.neurips.cc/paper_files/paper/2013/file/ac1dd209cbcc5e5d1c6e28598e8cbbe8-Paper.pdf}.

\bibitem[Kone{\v{c}}n{\`y} et~al.(2016)Kone{\v{c}}n{\`y}, McMahan, Yu,
  Richt{\'a}rik, Suresh, and Bacon]{konevcny2016federated}
J.~Kone{\v{c}}n{\`y}, H.~B. McMahan, F.~X. Yu, P.~Richt{\'a}rik, A.~T. Suresh,
  and D.~Bacon.
\newblock Federated learning: Strategies for improving communication
  efficiency.
\newblock \emph{arXiv preprint arXiv:1610.05492}, 2016.

\bibitem[Konečný et~al.(2017)Konečný, McMahan, Yu, Richtárik, Suresh, and
  Bacon]{konečný2017federated}
J.~Konečný, H.~B. McMahan, F.~X. Yu, P.~Richtárik, A.~T. Suresh, and
  D.~Bacon.
\newblock Federated learning: Strategies for improving communication
  efficiency, 2017.

\bibitem[Kovalev et~al.(2019)Kovalev, Horvath, and Richtarik]{kovalev2019dont}
D.~Kovalev, S.~Horvath, and P.~Richtarik.
\newblock Don't jump through hoops and remove those loops: Svrg and katyusha
  are better without the outer loop, 2019.

\bibitem[Li et~al.(2020)Li, Kovalev, Qian, and Richt{\'a}rik]{adiana}
Z.~Li, D.~Kovalev, X.~Qian, and P.~Richt{\'a}rik.
\newblock Acceleration for compressed gradient descent in distributed and
  federated optimization.
\newblock \emph{arXiv preprint arXiv:2002.11364}, 2020.

\bibitem[Li et~al.(2021)Li, Bao, Zhang, and Richtárik]{li2021page}
Z.~Li, H.~Bao, X.~Zhang, and P.~Richtárik.
\newblock Page: A simple and optimal probabilistic gradient estimator for
  nonconvex optimization, 2021.

\bibitem[Mishchenko et~al.(2023)Mishchenko, Gorbunov, Takáč, and
  Richtárik]{mishchenko2023distributed}
K.~Mishchenko, E.~Gorbunov, M.~Takáč, and P.~Richtárik.
\newblock Distributed learning with compressed gradient differences, 2023.

\bibitem[Nesterov(2012)]{doi:10.1137/100802001}
Y.~Nesterov.
\newblock Efficiency of coordinate descent methods on huge-scale optimization
  problems.
\newblock \emph{SIAM Journal on Optimization}, 22\penalty0 (2):\penalty0
  341--362, 2012.
\newblock \doi{10.1137/100802001}.
\newblock URL \url{https://doi.org/10.1137/100802001}.

\bibitem[Nguyen et~al.(2017)Nguyen, Liu, Scheinberg, and
  Tak{\'a}{\v{c}}]{nguyen2017sarah}
L.~M. Nguyen, J.~Liu, K.~Scheinberg, and M.~Tak{\'a}{\v{c}}.
\newblock Sarah: A novel method for machine learning problems using stochastic
  recursive gradient.
\newblock In \emph{International conference on machine learning}, pages
  2613--2621. PMLR, 2017.

\bibitem[Richtárik and Takáč(2013)]{richtárik2013distributed}
P.~Richtárik and M.~Takáč.
\newblock Distributed coordinate descent method for learning with big data,
  2013.

\bibitem[Smith et~al.(2018)Smith, Forte, Ma, Tak{\'a}{\v{c}}, Jordan, and
  Jaggi]{smith2018cocoa}
V.~Smith, S.~Forte, C.~Ma, M.~Tak{\'a}{\v{c}}, M.~I. Jordan, and M.~Jaggi.
\newblock Cocoa: A general framework for communication-efficient distributed
  optimization.
\newblock \emph{Journal of Machine Learning Research}, 18\penalty0
  (230):\penalty0 1--49, 2018.

\bibitem[Sun et~al.(2023)Sun, Xu, Yang, Nath, Li, Zhao, Xu, Chen, and
  Roth]{sun2023communicationefficient}
J.~Sun, Z.~Xu, D.~Yang, V.~Nath, W.~Li, C.~Zhao, D.~Xu, Y.~Chen, and H.~R.
  Roth.
\newblock Communication-efficient vertical federated learning with limited
  overlapping samples, 2023.

\bibitem[Szlendak et~al.(2021)Szlendak, Tyurin, and
  Richtárik]{szlendak2021permutation}
R.~Szlendak, A.~Tyurin, and P.~Richtárik.
\newblock Permutation compressors for provably faster distributed nonconvex
  optimization, 2021.

\bibitem[Verbraeken et~al.(2020)Verbraeken, Wolting, Katzy, Kloppenburg,
  Verbelen, and Rellermeyer]{verbraeken2020survey}
J.~Verbraeken, M.~Wolting, J.~Katzy, J.~Kloppenburg, T.~Verbelen, and J.~S.
  Rellermeyer.
\newblock A survey on distributed machine learning.
\newblock \emph{Acm computing surveys (csur)}, 53\penalty0 (2):\penalty0 1--33,
  2020.

\bibitem[Xu et~al.(2021)Xu, Fan, Li, and Yang]{xu2021efficient}
W.~Xu, H.~Fan, K.~Li, and K.~Yang.
\newblock Efficient batch homomorphic encryption for vertically federated
  xgboost, 2021.

\bibitem[Zhang et~al.(2021)Zhang, Gu, Deng, Gu, Bo, Pei, and Huang]{Zhang_2021}
Q.~Zhang, B.~Gu, C.~Deng, S.~Gu, L.~Bo, J.~Pei, and H.~Huang.
\newblock Asysqn: Faster vertical federated learning algorithms with better
  computation resource utilization.
\newblock In \emph{Proceedings of the 27th ACM SIGKDD Conference on Knowledge
  Discovery \& Data Mining}, pages 3917--3927, 2021.

\end{thebibliography}

%%%%%%%%%%%%%%%%%%%%%%%%%%%%%%%%%%%%%%%%%%%%%%%%%%%%%%%%%%%%
\appendix
\onecolumn
\newpage
\part*{Supplementary Material}
\tableofcontents

\newpage 

\section{Auxiliary Lemmas} \label{sec:aux}

\begin{lemma}
\textup{(Lemma 12 from \citep{kovalev2019dont})} For a random vector $x \in \mathbb{R}^d$ and any $y \in \mathbb{R}^d$, the variance of $y$ can be represented in a form:
\begin{equation}
    \label{lem:disp-prop}
    \mathbb{E} \left [\|x - \mathbb{E} \left [ x\right]\|^2 \right] = \mathbb{E} \left [ \|x - \|y\|^2 \right] - \mathbb{E} \left [ \| \mathbb{E} \left[ x\right ] - y \|^2\right]
\end{equation}
\end{lemma}
The next lemma is a form of Young’s inequality:
\begin{lemma}
 Let $\beta = \frac{\eta \theta_1}{L (1 - \eta \theta_1)}$. Then $\forall a, b \in \mathbb{R}^d$ the following inequality holds:
\begin{equation}
    \label{lem:young}
    \la a, b \ra \geq - \frac{\|a\|^2}{2 \beta} - \frac{\beta \|b\|^2}{2}
\end{equation}
\end{lemma}
% \setcitestyle{authoryear, round, comma, aysep={;}, yysep={,}, notesep={, }}
\section{Horizontal partition} \label{sec:horizontal}
Although this might seem trivial, we further rely on gradient unbiasedness. Thus,
\begin{lemma}\label{lem:horizontal-unbiasedness}
\textup{} In Algorithm~\ref{alg:one} $\mathbb{E}\left[ g^k \right] = \nabla f(x^k)$.
\end{lemma}
\textbf{Proof:}
From Algorithm~\ref{alg:one} we get:
\begin{equation}
\label{from-alg-1-get}
g^k \eqdef \frac{1}{n} \sum_{i = 1}^n Q_i \left ( \nabla f_{i} (x^k) - \nabla f_{i} (w^k) \right ) + \nabla f (w^k)
\end{equation}
Using Tower Property and Definition~\ref{ass:ass000} of compressors' unbiasedness:
\begin{equation*}
\mathbb{E} \left[g^k\right] \eqdef \frac{1}{n} \sum_{i = 1}^n \mathbb{E}\left[Q_i \left ( \nabla f_{i} (x^k) - \nabla f_{i} (w^k) \right )\right] + \nabla f (w^k) = \nabla f(x^k)
\end{equation*}
\begin{lemma}\label{lem:lemma6-horizontal}
\textup{(Revised Lemma 6 from \citep{kovalev2019dont}, proof for Lemma \ref{lemm:lemma2.4})} In Algorithm~\ref{alg:one} the following holds:
\begin{equation}
\label{lemma6-inequality}
\mathbb{E} \left[\| g^k - \nabla f(x^k) \|^2 \right] \leq \frac{2 L \omega}{n} \left ( f (w^k) - f(x^k) - \la \nabla f (x^k); w^k - x^k \ra \right )
\end{equation}
$$
\tilde L = \frac{L \omega}{n}
$$
\end{lemma}
\textbf{Proof:}
From Algorithm~\ref{alg:one}, we derive \ref{from-alg-1-get}, thus, with an empty term, after opening the square and taking advantage of the unbiased and independent compression operators in the scalar product:
\begin{eqnarray*}
\EE \left[ \left \| g^k - \nabla f(x^k) \right \|^2 \right]
&=& \mathbb{E} \left [ \left \| \frac{1}{n} \sum_{i = 1}^n Q^w_i \left \{ \nabla f_{i} (x^k) - \nabla f_{i} (w^k) \right \} + \underline{\nabla f (w^k)} - \underline{\nabla f(x^k)} \right \|^2 \right ]\\
&=& \EE \left[\left\| \frac{1}{n} \sum_{i = 1}^n \left( Q_i \left \{ \nabla f_{i} (x^k) - \nabla f_{i} (w^k) \right \} - \left \{ \underline{\nabla f_i (x^k)} - \underline{\nabla f_i (w^k)} \right \} \right) \right \|^2 \right]\\
&=& \frac{1}{n^2} \sum_{i = 1}^n \mathbb{E} \left [  \left \| Q_i \left \{ \nabla f_{i} (x^k) - \nabla f_{i} (w^k) \right \} - \left \{ \nabla f_i (x^k) - \nabla f_i (w^k) \right \}  \right \|^2 \right ]\\
&+& \frac{1}{n^2} \sum_{i \neq l} \mathbb{E} [ \la Q_i \left \{ \nabla f_{i} (x^k) - \nabla f_{i} (w^k) \right \} - \left \{ \nabla f_i (x^k) - \nabla f_i (w^k) \right \}, Q_l \{ \nabla f_{l} (x^k)\\
&-& \nabla f_{l} (w^k) \} - \left \{ \nabla f_l (x^k) - \nabla f_l (w^k) \right \} \ra ]\\
&=& \frac{1}{n^2} \sum_{i = 1}^n \mathbb{E} \left [ \left \| Q_i \left \{ \nabla f_{i} (x^k) - \nabla f_{i} (w^k) \right \} - \left \{ \nabla f_i (x^k) - \nabla f_i (w^k) \right \}  \right \|^2 \right ]
\end{eqnarray*}
Next, using Definition \ref{ass:ass000}:
\begin{eqnarray*}
\EE \left[ \left \| g^k - \nabla f(x^k) \right \|^2 \right]
&=& \frac{1}{n^2} \sum_{i = 1}^n \mathbb{E} \left [ \left \| Q_i \left \{ \nabla f_{i} (x^k) - \nabla f_{i} (w^k) \right \} - \left \{ \nabla f_i (x^k) - \nabla f_i (w^k) \right \}  \right \|^2 \right ]\\ 
&\leq& \frac{\omega}{n^2} \sum_{i = 1}^{n} \mathbb{E} \left [ \left \| \nabla f_{i} (x^k) - \nabla f_{i} (w^k) \right \|^2\right]
\end{eqnarray*}
And finally, using Lipschitz property (Assumption \ref{ass:ass001}) and convexity (Assumption \ref{ass:ass002}) of functions $f_i (x)$:
\begin{eqnarray*}
\EE \left[ \left \| g^k - \nabla f(x^k) \right \|^2 \right]
&\leq& \frac{\omega}{n^2} \sum_{i = 1}^n 2 L \left ( f_i (w^k) - f_i(x^k) - \la \nabla f_i (x^k); w^k - x^k \ra \right )\\
&\leq& \frac{2 L \omega}{n} \left ( f (w^k) - f(x^k) - \la \nabla f (x^k); w^k - x^k \ra \right )
\end{eqnarray*}
\EndProof{}

Lemmas 7, 8, as already stated in the main section, are obtained from the original lemmas by replacing the Lipschitz constant of the gradient with a constant $\frac{\omega}{n}$ times greater than the one given by Assumption \ref{ass:ass001}. We now need to show that nothing changes in how lemmas 7, 8 behave at their core, although an Algorithm has experienced slight changes:
\begin{lemma}\label{lem:lemma7}
\textup{(Revised Lemma 7 from \citep{kovalev2019dont})} In Algorithm~\ref{alg:one} the following holds:
\begin{equation}
\la g^k, x^* - z^{k+1} \ra + \frac{\mu}{2} \left \| x^k - x^* \right \|^2 \geq \frac{\widetilde{L}}{2 \mu} \left \| z^k - z^{k+1} \right \|^2 + Z^{k+1} - \frac{1}{1 + \eta \sigma} Z^k
\end{equation}
\end{lemma}
\textbf{Proof:}
From Algorithm \ref{alg:one}:
\begin{equation*}
z^{k+1} \eqdef \frac{1}{1 + \eta \frac{\mu}{\widetilde{L}}} \left (\eta \frac{\mu}{\widetilde{L}} x^k + z^k - \frac{\eta}{\widetilde{L}} g^k\right )
\end{equation*}
\begin{equation*}
\rightarrow g^k = \frac{\widetilde{L}}{\mu} \left (z^k - z^{k+1} \right ) + \mu \left (x^k - z^{k+1} \right )
\end{equation*}
Thus (we will simply neglect the underlined summand):
\begin{eqnarray*}
\la g^k, z^{k+1} - x^* \ra 
&=& \mu \la x^k - z^{k+1}, z^{k+1} - x^*\ra + \frac{\widetilde{L}}{\eta} \la z^k - z^{k+1}, z^{k+1} - x^*\ra\\
&=& \frac{\mu}{2} \left (\| x^k - x^* \|^2 - \underline{\| x^k - z^{k+1}  \|^2} - \|z^{k+1} - x^* \|^2  \right )\\
&+& \frac{\widetilde{L}}{2 \eta} \left (\| z^k - x^* \|^2 - \| z^k - z^{k+1} \|^2 - \| z^{k+1} - x^* \|^2 \right )\\
&\leq& \frac{\mu}{2} \| x^k - x^* \|^2 + \frac{\widetilde{L}}{2 \eta} \left (\|z^k - x^* \|^2 - (1 + \eta \sigma) \|z^{k+1} - x^*\|^2 \right ) - \frac{\widetilde{L}}{2 \eta} \|z^k - z^{k+1}\|^2
\end{eqnarray*}
\EndProof{Which concludes the proof }

Now let us show that Lemma 8 can be rewritten in the same manner, replacing $L$ - Lipschitz constant with the new $\widetilde{L}$.
\begin{lemma}\label{lem:lemma8}
\textup{(Revised Lemma 8 from \citep{kovalev2019dont}} We have:
\begin{equation}
\frac{1}{\theta_1} \left ( f(y^{k+1}) - f(x^k) \right ) - \frac{\theta_2}{2 \widetilde{L} \theta_1} \left \|g^k - \nabla f(x^k) \right \|^2 \leq \frac{\widetilde{L}}{2 \eta} \left \| z^{k+1} - z^k \right \|^2 + \la g^k , z^{k+1} - z^k \ra
\end{equation}
\end{lemma}

\textbf{Proof:}
Firstly, we will prove this lemma in the case of $\frac{\omega}{n} > 1$. From Algorithm \ref{alg:one}, one can see that $y^{k + 1} \gets x^k + \theta_1 (z^{k + 1} - z^k)$. Below we will introduce an empty term, use $\widetilde{L}$-smoothness of function $f(x)$ and, after that, utilize Lemma \ref{lem:young}, where $\beta = \frac{\eta \theta_1}{\widetilde{L} \left ( 1 - \eta \theta_1 \right ) }$:
\begin{eqnarray*}
\frac{\widetilde{L}}{2 \eta} \|z^{k+1} - z^k \|^2 + \la g^k, z^{k+1} - z^k\ra
&=& \frac{1}{\theta_1} \left (\frac{\widetilde{L}}{2 \eta \theta_1} \|\theta_1 \left (z^{k+1} - z^k \right ) \|^2 + \la g^k, \theta_1 \left (z^{k+1} - z^k \right )\ra \right )\\
&=& \frac{1}{\theta_1} \left (\frac{\widetilde{L}}{2 \eta \theta_1} \| y^{k+1} - x^k \|^2 + \la g^k, y^{k+1} - x^k\ra \right )\\
&=& \underline{\frac{\widetilde{L}}{2\theta_1} \| y^{k+1} - x^k \|^2} + \frac{1}{\theta_1} \la\underline{\nabla f(x^k)}, y^{k+1} - x^k\ra\\
&+&\frac{\widetilde{L}}{2\theta_1} \left (\frac{1}{\eta \theta_1} \underline{- 1} \right ) \| y^{k+1} - x^k \|^2 + \frac{1}{\theta_1} \la g^k \underline{- \nabla f(x^k)}, y^{k+1} - x^k\ra )\\
&\geq& \frac{1}{\theta_1} \left(f(y^{k+1}) - f(x^k)\right) + \frac{\widetilde{L}}{2\theta_1} \left (\frac{1}{\eta \theta_1} - 1 \right ) \| y^{k+1} - x^k \|^2 \\
&+& \frac{1}{\theta_1} \la g^k - \nabla f(x^k), y^{k+1} - x^k\ra\\
&\geq& \frac{1}{\theta_1} \left (f(y^{k+1}) - f(x^k) - \frac{\eta \theta_1}{2 \widetilde{L} ( 1 - \eta \theta_1)} \| g^k - \nabla f(x^k)\|^2 \right )\\
&=& \frac{1}{\theta_1} \left (f(y^{k+1}) - f(x^k) - \frac{\theta_2}{2 \widetilde{L}} \|g^k - \nabla f(x^k) \|^2\right )
\end{eqnarray*}

The last equality was received due to $\eta = \frac{\theta_2}{(1 + \theta_2) \theta_1}$. Thus, we derive:
\begin{equation*}
\frac{1}{\theta_1} \left (f(y^{k+1}) - f(x^k) \right ) - \frac{\theta_2}{2 \widetilde{L} \theta_1} \|g^k - \nabla f(x^k) \|^2 \leq \frac{\widetilde{L}}{2 \eta} \| z^{k+1} - z^k \|^2 + \la g^k, z^{k+1} - z^k \ra
\end{equation*}

This is a part of the lemma for the case $\frac{w}{n} \leq 1$, e.g. $\widetilde{L} \leq L$.
\begin{eqnarray*}
\frac{\widetilde{L}}{2 \eta} \|z^{k+1} - z^k \|^2 + \la g^k, z^{k+1} - z^k\ra
&=& \frac{1}{\theta_1} \left (\frac{\widetilde{L}}{2 \eta \theta_1} \|\theta_1 \left (z^{k+1} - z^k \right ) \|^2 + \la g^k, \theta_1 \left (z^{k+1} - z^k \right )\ra \right )\\
&=& \frac{1}{\theta_1} \left (\frac{\widetilde{L}}{2 \eta \theta_1} \| y^{k+1} - x^k \|^2 + \la g^k, y^{k+1} - x^k\ra \right )\\
&=&\frac{L}{2\theta_1} \| y^{k+1} - x^k \|^2 + \frac{1}{\theta_1} \la\nabla f(x^k), y^{k+1} - x^k\ra\\
&+&\frac{L}{2\theta_1} \left (\frac{1}{\frac{L}{\widetilde{L}}\eta \theta_1} - 1 \right ) \| y^{k+1} - x^k \|^2 + \frac{1}{\theta_1} \la g^k - \nabla f(x^k), y^{k+1} - x^k\ra )\\
&\geq& \frac{1}{\theta_1} \left(f(y^{k+1}) - f(x^k)\right) + \frac{L}{2\theta_1} \left (\frac{1}{\frac{L}{\widetilde{L}}\eta \theta_1} - 1 \right ) \| y^{k+1} - x^k \|^2 \\
&+& \frac{1}{\theta_1} \la g^k - \nabla f(x^k), y^{k+1} - x^k\ra\\
&\geq& \frac{1}{\theta_1} \left (f(y^{k+1}) - f(x^k) - \frac{\frac{L}{\widetilde{L}}\eta \theta_1}{2 \widetilde{L} ( 1 - \frac{L}{\widetilde{L}}\eta \theta_1)} \| g^k - \nabla f(x^k)\|^2 \right )\\
&=& \frac{1}{\theta_1} \left (f(y^{k+1}) - f(x^k) - \frac{\theta_2}{2 L} \|g^k - \nabla f(x^k) \|^2\right)
\end{eqnarray*}

The last inequality was obtained by using Young's inequality in the form of $\langle a, b \rangle \geq -\frac{\|a\|^2}{2\beta} - \frac{\beta\|b\|^2}{2}$ for $\beta = \frac{\frac{L}{\widetilde{L}}\eta \theta_1}{L\left(1 - \frac{L}{\widetilde{L}} \eta \theta_1\right)}$. \\
The last equality was received due to $\eta = \frac{\frac{\widetilde{L}}{L}\theta_2}{(1 + \theta_2) \theta_1}$. \\
Using that $\widetilde{L} \leq L$, we get that:
\begin{equation*}
\frac{1}{\theta_1} \left (f(y^{k+1}) - f(x^k) - \frac{\theta_2}{2 L} \|g^k - \nabla f(x^k) \|^2\right) \geq \frac{1}{\theta_1} \left (f(y^{k+1}) - f(x^k) - \frac{\theta_2}{2 \widetilde{L}} \|g^k - \nabla f(x^k) \|^2\right)
\end{equation*}
And finally we get:
\begin{equation*}
\frac{1}{\theta_1} \left (f(y^{k+1}) - f(x^k) \right ) - \frac{\theta_2}{2 \widetilde{L} \theta_1} \|g^k - \nabla f(x^k) \|^2 \leq \frac{\widetilde{L}}{2 \eta} \| z^{k+1} - z^k \|^2 + \la g^k, z^{k+1} - z^k \ra
\end{equation*}
\EndProof{}

Lemma 9, however, does not require any changes in the proof at all, because it is based only on the definition of $w^{k+1}$ from Algorithm \ref{alg:one}.
\begin{lemma}\label{lem:lemma9}
\textup{(Lemma 9 from \citep{kovalev2019dont})} We have:
\begin{equation}
\mathbb{E} \left [ f(w^{k+1} \right] = \left (1 - p \right) f(w^k) + p f(y^k)
\end{equation}
\end{lemma}
\begin{lemma}\label{lem:lemma10-horizontal}
\textup{(Revised Lemma 10 from \citep{kovalev2019dont}, proof for Theorem \ref{th:one})} Considering lemmas \ref{lem:lemma6-horizontal}–\ref{lem:lemma9}, we get:
\begin{equation}
\mathbb{Z}^k\left[\frac{1}{1 + \eta\sigma} \right] + \mathbb{Y}^k\left[(1 - \theta_1(1 - \theta_2)) \right] + \mathbb{W}^k\left[1 - \frac{p\theta_1}{1 + \theta_1}\right] \geq \mathbb{E}\left[\mathbb{Z}^{k + 1} + \mathbb{Y}^{k + 1} + \mathbb{W}^{k + 1} \right]
\end{equation}
\end{lemma}
\textbf{Proof:} Using strong convexity (Assumption \ref{ass:ass002}) and adding an empty term, we get:
\begin{eqnarray*}
f(x^*)
&\geq& f(x^k) + \la\nabla f(x^k) , x^* - x^k\ra + \frac{\mu}{2} \|x^k - x^* \|^2\\
&\geq& f(x^k) + \la\nabla f(x^k) , x^* \underline{- z^k} + \underline{z^k} - x^k\ra + \frac{\mu}{2} \|x^k - x^* \|^2
\end{eqnarray*}

From Algorithm \ref{alg:one}, $z^k - x^k \eqdef \frac{\theta_2}{\theta_1} \left ( x^k - w^k \right ) + \frac{1 - \theta_1 - \theta_2}{\theta_1} \left ( x^k - y^k \right ) $, thus:
\begin{eqnarray*}
f(x^*)
&\geq& f(x^k) + \la\nabla f(x^k) , x^* \underline{- z^k} + \underline{z^k} - x^k\ra + \frac{\mu}{2} \|x^k - x^* \|^2\\
&\geq& f(x^k) + \frac{\mu}{2} \|x^k - x^* \|^2 + \la \nabla f(x^k) , x^* - z^k \ra\\
&+& \frac{\theta_2}{\theta_1} \la \nabla f(x^k) , x^k - w^k \ra + \underline{\frac{1 - \theta_1 - \theta_2}{\theta_1} \la \nabla f(x^k) , x^k - y^k\ra}
\end{eqnarray*}

Applying convexity to the underlined term, utilizing unbiasedness of the gradient and adding an empty term, we derive:
\begin{eqnarray*}
f(x^*)
&\geq& f(x^k) + \frac{\mu}{2} \|x^k - x^* \|^2 + \la \nabla f(x^k) , x^* - z^k \ra\\
&+& \frac{\theta_2}{\theta_1} \la \nabla f(x^k) , x^k - w^k \ra + \underline{\frac{1 - \theta_1 - \theta_2}{\theta_1} \la \nabla f(x^k) , x^k - y^k\ra}\\
&\geq& f(x^k) + \frac{\theta_2}{\theta_1} \la \nabla f(x^k) , x^k - w^k \ra + \frac{1 - \theta_1 - \theta_2}{\theta_1} \left ( f(x^k) - f(y^k) \right )\\
&+& \mathbb{E} \left [ \frac{\mu}{2} \|x^k - x^*\|^2 \la g^k, x^* \underline{- z^{k+1}}\ra + \la g^k, \underline{z^{k+1}} - z^k\ra \right ]\\
&=& f(x^k) + \frac{\theta_2}{\theta_1} \la \nabla f(x^k) , x^k - w^k \ra + \frac{1 - \theta_1 - \theta_2}{\theta_1} \left ( f(x^k) - f(y^k) \right )\\
&+& \mathbb{E} \left [\frac{\mu}{2} \|x^k - x^*\|^2 + \la g^k, x^* - z^{k+1}\ra + \la, z^{k+1} - z^k\ra  \right ]
\end{eqnarray*}

Using Lemma \ref{lem:lemma7}, we get:
\begin{eqnarray*}
f(x^*)
&=& f(x^k) + \frac{\theta_2}{\theta_1} \la \nabla f(x^k) , x^k - w^k \ra + \frac{1 - \theta_1 - \theta_2}{\theta_1} \left ( f(x^k) - f(y^k) \right )\\
&+& \mathbb{E} \left [\frac{\mu}{2} \|x^k - x^*\|^2 + \la g^k, x^* - z^{k+1}\ra + \la g^k, z^{k+1} - z^k\ra  \right ]\\
&\geq& f(x^k) + \frac{\theta_2}{\theta_1} \la \nabla f(x^k) , x^k - w^k \ra + \frac{1 - \theta_1 - \theta_2}{\theta_1} \left ( f(x^k) - f(y^k) \right )\\
&+& \mathbb{E} \left [ Z^{k+1} - \frac{1}{1 + \eta \sigma} Z^k \right] + \mathbb{E} \left [ \la g^k , z^{k+1} - z^k \ra + \frac{\widetilde{L}}{2 \eta} \|z^k - z^{k+1} \|^2 \right]\\
&=& f(x^k) + \frac{\theta_2}{\theta_1} \la \nabla f(x^k) , x^k - w^k \ra + \frac{1 - \theta_1 - \theta_2}{\theta_1} \left ( f(x^k) - f(y^k) \right )\\
&+& \mathbb{E} \left [ Z^{k+1} - \frac{1}{1 + \eta \sigma} Z^k \right] + \mathbb{E} \underline{\left [ \la g^k , z^{k+1} - z^k \ra + \frac{\widetilde{L}}{2 \eta} \|z^k - z^{k+1} \|^2 \right]}
\end{eqnarray*}

Utilizing Lemma \ref{lem:lemma8}, we transform an underlined summand:
\begin{eqnarray*}
f(x^*)
&=& f(x^k) + \frac{\theta_2}{\theta_1} \la \nabla f(x^k) , x^k - w^k \ra + \frac{1 - \theta_1 - \theta_2}{\theta_1} \left ( f(x^k) - f(y^k) \right )\\
&+& \mathbb{E} \left [ Z^{k+1} - \frac{1}{1 + \eta \sigma} Z^k \right] + \mathbb{E} \underline{\left [ \la g^k , z^{k+1} - z^k \ra + \frac{\widetilde{L}}{2 \eta} \|z^k - z^{k+1} \|^2 \right]}\\
&\geq& f(x^k) + \frac{\theta_2}{\theta_1} \la \nabla f(x^k) , x^k - w^k \ra + \frac{1 - \theta_1 - \theta_2}{\theta_1} \left ( f(x^k) - f(y^k) \right )\\
&+& \mathbb{E} \left [ Z^{k+1} - \frac{1}{1 + \eta \sigma} Z^k \right] + \mathbb{E} \underline{\left [ \frac{1}{\theta_1} \left ( f(y^{k+1}) - f(x^k) \right ) - \frac{\theta_2}{2 \widetilde{L} \theta_1} \|g^k - \nabla f(x^k) \|^2 \right ]}\\
&\geq& f(x^k) + \frac{\theta_2}{\theta_1} \la \nabla f(x^k) , x^k - w^k \ra + \frac{1 - \theta_1 - \theta_2}{\theta_1} \left ( f(x^k) - f(y^k) \right )\\
&+& \mathbb{E} \left [ Z^{k+1} - \frac{1}{1 + \eta \sigma} Z^k + \frac{1}{\theta_1} \left ( f(y^{k+1}) - f(x^k) \right ) \right] - \mathbb{E} \left [ 
 \frac{\theta_2}{2 \widetilde{L} \theta_1} \|g^k - \nabla f(x^k) \|^2 \right ]
\end{eqnarray*}

Finally, using Lemma \ref{lem:lemma6-horizontal} to evaluate the last summ and, we get:
\begin{eqnarray*}
f(x^*)
&\geq& f(x^k) + \frac{\theta_2}{\theta_1} \la \nabla f(x^k) , x^k - w^k \ra + \frac{1 - \theta_1 - \theta_2}{\theta_1} \left ( f(x^k) - f(y^k) \right )\\
&+& \mathbb{E} \left [ Z^{k+1} - \frac{1}{1 + \eta \sigma} Z^k + \frac{1}{\theta_1} \left ( f(y^{k+1}) - f(x^k) \right ) \right] - \mathbb{E} \left [ 
 \frac{\theta_2}{2 \widetilde{L} \theta_1} \|g^k - \nabla f(x^k) \|^2 \right ]\\ 
&\geq& f(x^k) + \frac{\theta_2}{\theta_1} \la \nabla f(x^k) , x^k - w^k \ra + \frac{1 - \theta_1 - \theta_2}{\theta_1} \left ( f(x^k) - f(y^k) \right )\\
&+& \mathbb{E} \left [ Z^{k+1} - \frac{1}{1 + \eta \sigma} Z^k + \frac{1}{\theta_1} \left ( f(y^{k+1}) - f(x^k) \right ) \right]\\
&-&\frac{\theta_2}{2 \widetilde{L} \theta_1} \mathbb{E} \underline{\left [ 
 2 \widetilde{L} \left ( f (w^k) - f(x^k) - \la \nabla f (x^k); w^k - x^k \ra \right ) \right ]}\\
&=& f(x^k) + \frac{\theta_2}{\theta_1} \la \nabla f(x^k) , x^k - w^k \ra + \frac{1 - \theta_1 - \theta_2}{\theta_1} \left ( f(x^k) - f(y^k) \right )\\
&+& \mathbb{E} \left [ Z^{k+1} - \frac{1}{1 + \eta \sigma} Z^k + \frac{1}{\theta_1} \left ( f(y^{k+1}) - f(x^k) \right ) \right]\\
&-&\frac{\theta_2}{\theta_1} \mathbb{E} \left [ f (w^k) - f(x^k) - \la \nabla f (x^k); w^k - x^k \ra \right ]\\
&=& - \frac{\theta_2}{\theta_1} f(w^k) - \frac{1 - \theta_1 - \theta_2}{\theta_1} f(y^k) + \left [ Z^{k+1} - \frac{1}{1 + \eta \sigma} Z^k + \frac{1}{\theta_1} f(y^{k+1}) \right]
\end{eqnarray*}

As one can see, this result corresponds exactly to the formulation of Lemma 10 from the original article, with the only difference that now everywhere we have the constant $\widetilde{L}$ instead of $L$, which leads us to the same final expression:
\begin{equation}
\mathbb{Z}^k\left[\frac{1}{1 + \eta\sigma} \right] + \mathbb{Y}^k\left[(1 - \theta_1(1 - \theta_2)) \right] + \mathbb{W}^k\left[1 - \frac{p\theta_1}{1 + \theta_1}\right] \geq \mathbb{E}\left[\mathbb{Z}^{k + 1} + \mathbb{Y}^{k + 1} + \mathbb{W}^{k + 1} \right]
\end{equation}\EndProof{}

From Lemma \ref{lem:lemma10-horizontal} we get our main Theorem directly by substituting into the inequality the specified values that preserve the inequality:
\begin{theorem}\label{theorem:theorem-horizontal}
\textup{(Theorem 11 from \citep{kovalev2019dont})} Let Assumptions~\ref{ass:ass001},~\ref{ass:ass002} be hold. Additionally, let $p = \frac{1}{n}$, $\theta_1 = \min{\{ \sqrt{\frac{2\sigma n }{3}} , \frac{1}{2} \}}$, $\theta_2 = \frac{1}{2}$. Then $\mathbb{E}\left[\Psi^k\right] \leq \epsilon \Psi^0$ after $K = \mathcal{O}\left(\left(n + \sqrt{\frac{\omega L}{\mu}}\right)\log{\frac{1}{\epsilon}}\right)$ iterations.
\end{theorem}

\begin{theorem}\label{theorem:PermK horizontal proof}
\textup{(The efficient Lipschitz for PermK compressor in the horizontal case, proof for Lemma \ref{lem:PermK horizontal})}\\
For \texttt{DHPL-Katyusha} with PermK compressor, the efficient  Lipschitz constant is equal to $L$.
\end{theorem}
\textbf{Proof:}
Denote $Q_i$ as a PermK compressor that operates on $i$-th function. Therefore, we get:
\begin{eqnarray*}
\EE \left[ \left \| g^k - \nabla f(x^k) \right \|^2 \right]
&=& \mathbb{E} \left [ \left \| \frac{1}{n} \sum_{i = 1}^n Q^w_i \left \{ \nabla f_{i} (x^k) - \nabla f_{i} (w^k) \right \} + \nabla f (w^k) - \nabla f(x^k) \right \|^2 \right ]\\
&=& \EE \left[\left\| \frac{1}{n} \sum_{i = 1}^n  Q_i \left \{ \nabla f_{i} (x^k) - \nabla f_{i} (w^k) \right \}  - \frac{1}{n}\sum_{i = 1}^n\left \{ \nabla f_i (x^k) - \nabla f_i (w^k) \right \} \right \|^2 \right]\\
&\leq& \EE \left [A\frac{1}{n} \sum_{i = 1}^n \left \|\nabla f_i (x^k) -\nabla f_i (w^k) \right \|^2 - B\left \|  \frac{1}{n}\sum_{i = 1}^n\nabla f_i (x^k) - \nabla f_i (w^k) \right \|^2 \right]\\
&\leq&  \EE \left[A\frac{1}{n} \sum_{i = 1}^n \left \| \nabla f_i (x^k) - \nabla f_i (w^k) \right \|^2 \right] \\
&\leq& 2LA \left ( f (w^k) - f(x^k) - \la \nabla f (x^k), w^k - x^k \ra \right )
\end{eqnarray*}
The first inequality is obtained by using A-B inequality from \citep{szlendak2021permutation}, the second inequality we get by utilizing the positive definiteness of the norm and the last inequality is obtained by using an analogy with \ref{lem:horizontal-unbiasedness}.

In the case $n \geq d$ we have $A = \frac{d - 1}{n - 1} \leq 1$, and in the case $n < d$ we have $A = 1$. Combining these cases, we get the estimation for the efficient Lipschitz constant as $1$.

\newpage
\section{Vertical case}\label{sec:vertical}
In this section, we provide the complete proofs for our results in 
the vertical case.

\begin{lemma}
\textup{(Revised Lemma 6 from \citep{kovalev2019dont}, proof for Lemma \ref{lem: lemma 6 arbitrary loss})}\label{lemma6 arbitrary vertical} In Algorithm~\ref{alg:orig} the following holds:
\begin{eqnarray*}
&\mathbb{E}\left\|g^k(x^k) - \nabla\mathcal{L}\left(Ax^k, b_j\right)\right\|^2 \\
&\leq 2\max \left\{L, \frac{\bar{L}}{K} \right\} \left(\nabla\mathcal{L}\left(Aw^k, b\right) - \nabla\mathcal{L}\left(Ax^k, b\right) - \langle\nabla\mathcal{L}\left(Ax^k, b\right); w^k - x^k \rangle\right)
\end{eqnarray*}
\end{lemma}
\textbf{Proof:}
\begin{eqnarray*}
\mathbb{E}||g^k(x^k) &-& \nabla\mathcal{L}\left(Ax^k, b_j\right)||^2 \\
&=& \mathbb{E}_{J^k}\left\| \frac{1}{K}\sum\limits_{j \in {J^k}}\frac{1}{sp_j}\left[\nabla\mathcal{L}_j\left(\sum\limits_{i = 1}^{n}X_{ij}, b_j\right)\!-\!\nabla\mathcal{L}_j\left(\sum\limits_{i = 1}^{n}W_{ij}, b_j\right)\right] + \nabla\mathcal{L}\left(Aw^k, b_j\right)\!-\!\nabla\mathcal{L}\left(Ax^k, b_j\right)\right\|^2\\
&=&\frac{1}{K}\!\mathbb{E}_{i \sim D}\!\left\|\!\left(\frac{1}{sp_j}\!\left[\nabla\mathcal{L}_j\left(\sum\limits_{i = 1}^{n}X_{ij}, b_j\right)\!-\!\nabla\mathcal{L}_j\left(\sum\limits_{i = 1}^{n}W_{ij}, b_j\right)\right]\right)\!-\! \left(\nabla\mathcal{L}\left(Aw^k, b_j\right)\!+\!\nabla\mathcal{L}\left(Ax^k, b_j\right)\right)\right\|^2 \\
&\leq&\frac{1}{K}\mathbb{E}_{i \sim D}\left\| \frac{1}{sp_j}\left[\nabla\mathcal{L}_j\left(\sum\limits_{i = 1}^{n}X_{ij}, b_j\right)-\nabla\mathcal{L}_j\left(\sum\limits_{i = 1}^{n}W_{ij}, b_j\right)\right]\right\|^2\\
&\leq&\frac{1}{K}\sum\limits_{j = 1}^s \frac{2L_j}{s^2 p_j} \left( \nabla\mathcal{L}\left(Aw^k, b\right) - \nabla\mathcal{L}\left(Ax^k, b\right) - \langle\nabla\mathcal{L}\left(Ax^k, b\right); w^k - x^k \rangle\right) \\
&=& \frac{2\bar{L}}{K}\left( \nabla\mathcal{L}\left(Aw^k, b\right) - \nabla\mathcal{L}\left(Ax^k, b\right) - \langle\nabla\mathcal{L}\left(Ax^k, b\right); w^k - x^k \rangle\right) \\
&\leq& 2\max \left\{L, \frac{\bar{L}}{K} \right\} \left( \nabla\mathcal{L}\left(Aw^k, b\right) - \nabla\mathcal{L}\left(Ax^k, b\right) - \langle\nabla\mathcal{L}\left(Ax^k, b\right); w^k - x^k \rangle\right)
\end{eqnarray*}

\begin{lemma}
\textup{(Revised Lemma 6 from \citep{kovalev2019dont}, proof for Lemma \ref{lem: lemma 6 vertically distributed})}\label{lemma6vert} In Algorithm~\ref{alg:two} the following holds:
\begin{eqnarray*}
\mathbb{E}||g^k(x^k) - \nabla f(x^k)||^2 &\leq& 2L \left( f (w^k) - f(x^k) - \langle \nabla f (x^k); w^k - x^k \rangle\right)\\
&\cdot& \left(1 + (\omega - 1)\frac{s \sum\limits_{j = 1}^s L_j^2}{\mu^2} \right)
\end{eqnarray*}
\end{lemma}
\textbf{Proof:}
\begin{eqnarray*}
\mathbb{E} \left [ \left \|g^k(x^k) - \nabla f(x^k) \right \|^2 \right]
&=& \mathbb{E}\left[ \left \| \frac{1}{b_s} \sum_{j \in \{J\}} \left [ 2 \sum_{i = 1}^n Q_i\left(\langle A^T_{ji}, x^k_i - w^k_i \rangle \right) \right ]A^T_j + \nabla f(w^k) - \nabla f(x^k)\right \|^2\right]\\
&\leq& \frac{1}{b_s^2} \cdot b_s \sum\limits_{j \in J} \cdot \mathbb{E}\left[ \left \|2 \sum_{i = 1}^n Q_i\left(\langle A^T_{ji}, x^k_i - w^k_i \rangle \right) A^T_j + \nabla f(w^k) - \nabla f(x^k) \right \|^2\right]
\end{eqnarray*}
Here we used a Cauchy-Bunyakovsky-Schwarz inequality. The next step is to evaluate the term within mathematic expectancy.
\begin{eqnarray*}
&&\mathbb{E}\left[ \left \|2 \sum_{i = 1}^n Q_i\left(\langle A^T_{ji}, x^k_i - w^k_i \rangle \right) A^T_j + \nabla f(w^k) - \nabla f(x^k) \right\|^2\right]\\ 
&=& \mathbb{E}  \left [ \left \| 2 A_j^T \sum_{i = 1}^n Q_i\left(\langle A^T_{ji}, x^k_i - w^k_i \rangle \right) - \frac{2}{s} A^T A(x^k - w^k) \right \|^2 \right ]\\
&=& \mathbb{E} \left [ \left \|2A_j^T \sum_{i = 1}^n Q_i\left(\langle A^T_{ji}, x^k_i - w^k_i \rangle \right) \right \|^2 \right] + \mathbb{E} \left [ \left \|\frac{2}{s} A^T A (x^k - w^k)\right \|^2 \right ]\\
&-& 2 \cdot \mathbb{E} \left [ \langle 2 A_j^T  \sum_{i = 1}^n Q_i\left(\langle A^T_{ji}, x^k_i - w^k_i \rangle \right), \frac{2}{s} A^T A (x^k - w^k) \rangle \right]\\ 
&=& \mathbb{E} \left [\left \|2 A_j^T \sum_{i = 1}^n Q_i\left(\langle A^T_{ji}, x^k_i - w^k_i \rangle \right) \right \|^2 \right ] - \left \|\frac{2}{s} A^T A (x^k - w^k)\right \|^2\\
&\leq& \mathbb{E} \left [\left \|2 A_j^T \sum_{i = 1}^n Q_i\left(\langle A^T_{ji}, x^k_i - w^k_i \rangle \right) \right \|^2 \right ]
\end{eqnarray*}

Now, using norm properties and Definition \ref{ass:ass000} we get:
\begin{eqnarray*}
&&\mathbb{E} \left [\left \|2 A_j^T \sum_{i = 1}^n Q_i\left(\langle A^T_{ji}, x^k_i - w^k_i \rangle \right) \right \|^2 \right ]\\
&=& 4 \mathbb{E} \left [\left \|A_j^T \right\|^2 \left \vert \sum_{i = 1}^n Q_i\left(\langle A^T_{ji}, x^k_i - w^k_i \rangle \right) \right \vert^2 \right ] \\
&=& 4 \mathbb{E} \left [\left \|A_j^T \right \|^2 \left \vert \sum_{i = 1}^n \left(Q_i\left(\langle A^T_{ji}, x^k_i\!-\!w^k_i \rangle \right)\right)^2 + \sum\limits_{i \neq t} Q_i\left(\langle A^T_{ji}, x^k_i\!-\!w^k_i \rangle \right) Q_t\left(\langle A^T_{jt}, x^k_t\!-\!w^k_t \rangle \right) \right \vert \right ]\\
&\leq& 4 \mathbb{E} \left [\left \| A_j^T \right \|^2 \left \vert \sum_{i = 1}^n \omega (\langle A^T_{ji}, x^k_i - w^k_i \rangle)^2  + \sum\limits_{i \neq t} \left(\langle A^T_{ji}, x^k_i - w^k_i \rangle \right) \left(\langle A^T_{jt}, x^k_t - w^k_t \rangle \right) \right \vert \right ] \\
&=& 4 \mathbb{E} \left[ \left\|A_j^T \right \|^2 \left \vert \sum_{i = 1}^n (\omega - 1) (\langle A^T_{ji}, x^k_i - w^k_i \rangle)^2  + \left(\langle A^T_j, x^k - w^k \rangle \right)^2 \right \vert \right] \\
&=& \mathbb{E} \left[ \left \|\nabla f_j(x^k) - \nabla f_j(w^k) \right \|^2 + 4(\omega - 1) \left\|A_j^T \right \|^2 \sum\limits_{i = 1}^n \left \vert\langle A^T_{ji}, x^k_i - w^k_i \rangle \right \vert^2\right]\\
&\leq& \mathbb{E} \left[ \left \|\nabla f_j(x^k) - \nabla f_j(w^k) \right \|^2 + 4(\omega - 1) \left \|A_j^T \right \|^4 \left \| x^k - w^k \right \|^2\right] \\
&=& \mathbb{E} \left[ \left \|\nabla f_j(x^k) - \nabla f_j(w^k) \right \|^2 + 4(\omega - 1) \left\|A_j^T \right \|^4 \left \| \left(\frac{2}{s}A^TA\right)^{-1}(\nabla f(x^k) - \nabla f(w^k)) \right \|^2\right]\\
&\leq& \mathbb{E}\left[ \left \|\nabla f_j(x^k) - \nabla f_j(w^k) \right \|^2 + (\omega - 1) \left\|A_j^T \right \|^4 s^2 \left \| \left(A^TA\right)^{-1} \right \|^2 \left \| \nabla f(x^k) - \nabla f(w^k) \right \|^2\right]\\
&\leq& \mathbb{E}\left[ \left \|\nabla f_j(x^k) - \nabla f_j(w^k) \right \|^2 + (\omega - 1) \left\|A_j^T \right \|^4 \frac{s^2}{\mu^2} \left \| \nabla f(x^k) - \nabla f(w^k) \right \|^2\right]\\
&\leq& \left ( f (w^k) - f(x^k) - \la \nabla f (x^k); w^k - x^k \ra \right ) \left(2L + 2L(\omega - 1)\frac{s^2 \cdot \frac{\sum\limits_{j = 1}^s L_j^2}{s}}{\mu^2} \right)\\
&=& 2L\left( f (w^k) - f(x^k) - \la \nabla f (x^k); w^k - x^k \ra \right ) \left(1 + (\omega - 1)\frac{s \cdot \sum\limits_{j = 1}^s L_j^2}{\mu^2} \right)\\
\end{eqnarray*}
Where the first equality is obtained by using absolute homogeneity of vector norm and having that $\sum_{i = 1}^n Q_i\left(\langle A^T_{ji}, x^k_i - w^k_i \rangle \right)$ is a scalar value, the first inequality is obtained by using compressor property \ref{ass:ass000},  the second and the third inequalities are obtained by using the definition of the norm as the supremum,  the fourth inequality is obtained by using that the spectral norm  of matrix is bounded by the biggest eigenvalue, which is bounded by $\frac{1}{\mu}$, final inequality is obtained by using the Lipschitz property \ref{ass:ass001}.
\\

\begin{lemma} \textup{(Proof for Lemma \ref{lem: lower bound vertical})} There exists such function that holds under \ref{ass:ass005} that its efficient Lipschitz constant $\widetilde{L}$ is proportional to  $(\omega - 1) \frac{L^3}{\mu^2}$.
\label{lem: lemma 6 vertically distributed proof} 
\end{lemma}
\textbf{Proof:}\\
Assume the size of batch $b = 1$. Therefore, in previous lemma we have an exact quality:
\begin{eqnarray*}
\mathbb{E} \left [\|g^k - \nabla f(x^k) \|^2 \right] = 4 \mathbb{E} \left[ \left \|\nabla f_j(x^k) - \nabla f_j(w^k) \right \|^2 + (\omega - 1) \left\|A_j^T \right \|^2 \sum\limits_{i = 1}^n \left |\langle a^i_j, x^k_i - w^k_i \rangle \right |^2\right]
\end{eqnarray*}
Now assume matrix $A$ is given as ($b < a$):
\begin{center}
$A = \begin{pmatrix}
a & a\\
b & -b
\end{pmatrix}$
\end{center}
Assume that $k = 2$ and on the first iteration of \texttt{DVPL-Katyusha} $\nabla f_2$ was chosen and $b$ (from $Ax - b$) and $x_0$ are collinear with $A_2^T$.\\
Denote $x^k - w^k = (c, -c)^T$. It has its first coordinate equal to the negative second, because it is orthogonal to $A_1^T$.\\

As $A_1^T$ is orthogonal with $x^k - w^k$:
$$\nabla f_1 (x^k) - \nabla f_1(w^k) = A_1^T \langle A_1^T, x^k - w^k \rangle = 0$$\\

Therefore, we have the difference of the target gradient in $x^k$ and in $w^k$ to be equal to:
$$\nabla f (x^k) - \nabla f(w^k) = \nabla f_1 (x^k) - \nabla f_1(w^k) + \nabla f_2 (x^k) - \nabla f_2(w^k) = \nabla f_2 (x^k) - \nabla f_2(w^k)$$\\

Note that the norm of $w^k - x^k$ equals to:
$$||w^k - x^k||^2 = 2c^2 = \frac{||\nabla f_2(x^k) - \nabla f_2(w^k)||^2}{4||A_2^T||^4} = \frac{||\nabla f(x^k) - \nabla f(w^k)||^2}{4||A_2^T||^4}$$

Also note that eigenvalues of $A^TA$ are equal to $2a^2$ and $2 b^2$. Therefore $L = 2a^2$ and $\mu = 2b^2$.

Combining all of these facts we get:

\begin{eqnarray*}
&&\mathbb{E} \left [\|g^k - \nabla f(x^k) \|^2 \right] \\
&=& 4 \mathbb{E} \left[ \left \|\nabla f_j(x^k) - \nabla f_j(w^k) \right \|^2 + (\omega - 1) \left\|A_j^T \right \|^2 \sum\limits_{i = 1}^n \left |\langle a^i_j, x^k_i - w^k_i \rangle \right |^2\right]\\
&=& 4 \mathbb{E} \left[ \left \|\nabla f_j(x^k) - \nabla f_j(w^k) \right \|^2 \right] + (\omega - 1) \left\|A_1^T \right \|^2 \left(a^2c^2 + a^2c^2 \right) + (\omega - 1) \left\|A_2^T \right \|^2 \left(b^2c^2 + b^2c^2 \right)\\
&=&  4 \mathbb{E} \left[ \left \|\nabla f_j(x^k) - \nabla f_j(w^k) \right \|^2 \right] + 2(\omega - 1) \left\|A_1^T \right \|^4 c^2 + 2(\omega - 1) \left\|A_2^T \right \|^4 c^2 \\
&=& 4 \mathbb{E} \left[ \left \|\nabla f_j(x^k) - \nabla f_j(w^k) \right \|^2 \right] + 2(\omega - 1) \left\|A_1^T \right \|^4 c^2 + 2(\omega - 1) \left\|A_2^T \right \|^4 c^2\\
&=& 4 \mathbb{E} \left[ \left \|\nabla f_j(x^k) - \nabla f_j(w^k) \right \|^2 \right] + (\omega - 1) \frac{L^2}{2\mu^2} ||\nabla f(x^k) - \nabla f(w^k)||^2 + \frac{1}{2}(\omega - 1)||\nabla f(x^k) - \nabla f(w^k)||^2
\end{eqnarray*}
Using the Lipschitz property in analogue with \ref{lemma6vert}, we conclude the proof.

 \subsection{PermK compressor}
\begin{lemma} \textup{(Proof for Lemma \ref{lem: vertical PermK})}
The efficient Lipschitz constant for PermK compressor is less than:
\begin{align*}
\widetilde{L}_{Perm} = {2 L} \frac{s\sum\limits_{j = 1}^s L_j^2}{\mu^2}.
\end{align*}
\end{lemma}
\textbf{Proof:}
\begin{eqnarray*}
\mathbb{E}\left\|g^k - \nabla f(x^k)\right\|^2 &=& \mathbb{E}\left\|2\frac{n}{n}\sum\limits_{i = 1}^n A_i^T \sum\limits_{j = 1}^n \langle A^T_{ji}, x^k_{ji} - w^k_{ji} \rangle I_{ij} - \frac{2}{s}A^TA\left(x^k - w^k\right)\right\|^2\\
&=& 4\mathbb{E}\left\|\sum\limits_{i = 1}^n A_i^T \sum\limits_{j = 1}^n \langle A^T_{ji}, x^k_{ji} - w^k_{ji} \rangle I_{ij} - \frac{1}{s}A^TA\left(x^k - w^k\right)\right\|^2\\
&=& 4\mathbb{E}\left\|\sum\limits_{i = 1}^n \left[A_i^T \sum\limits_{j = 1}^n \langle A^T_{ji}, x^k_{ji} - w^k_{ji} \rangle I_{ij} - \frac{1}{ns}A^TA\left(x^k - w^k\right)\right]\right\|^2 \\
&\leq& 8\sum\limits_{i = 1}^n\mathbb{E}\left\|A_i^T \sum\limits_{j = 1}^n \langle A^T_{ji}, x^k_{ji} - w^k_{ji} \rangle I_{ij} - \frac{1}{ns}A^TA\left(x^k - w^k\right)\right\|^2\\
&=& 8\sum\limits_{i = 1}^n\mathbb{E}\left[\left\|A_i^T \sum\limits_{j = 1}^n\langle A^T_{ji}, x^k_{ji} - w^k_{ji} \rangle I_{ij}\right\|^2 - \left\|\frac{1}{ns}A^TA\left(x^k - w^k\right)\right\|^2\right] \\
&\leq&  8\sum\limits_{i = 1}^n\mathbb{E}\left[\left(\sum\limits_{j = 1}^n \langle A^T_{ji}, x^k_{ji} - w^k_{ji} \rangle I_{ij} \right)^2 \left\|A_i^T\right\|^2\right]\\
&=&  8\sum\limits_{i = 1}^n\mathbb{E}\left[\left(\sum\limits_{j = 1}^n (\langle A^T_{ji},\!x^k_{ji}\!-\!w^k_{ji} \rangle I_{ij})^2\!+\!\sum\limits_{j \neq t} (\langle A^T_{ji}, x^k_{ji}\!-\!w^k_{ji} \rangle I_{ij})(\langle A^T_{jt}, x^k_{ti}\!-\!w^k_{ti} \rangle I_{it})  \right) \cdot\left\|A_i^T\right\|^2\right] =\\
\end{eqnarray*}
Using that $\mathbb{E} I^2_{ij} = \frac{1}{n}$ and as any sample can be chosen only by one worker, we get that $\mathbb{E} I_{it} I_{ij} = 0, j \neq t$. Therefore, using the analogy with \ref{lemma6vert} we get:

\begin{eqnarray*}
&=&  8\sum\limits_{i = 1}^n\mathbb{E}\left[\frac{1}{n}\sum\limits_{j = 1}^n \langle A^T_{ji}, x^k_{ji} - w^k_{ji} \rangle^2 \left\|A_i^T\right\|^2\right]\\
&\leq& {4 L} \left ( f (w^k) - f(x^k) - \langle \nabla f (x^k); w^k - x^k \rangle \right )\frac{s\sum\limits_{j = 1}^s L_j^2}{\mu^2}
\end{eqnarray*}

\end{document}